\documentclass{amsart}
\parskip=\smallskipamount

\newtheorem{theorem}{Theorem}[section]
\newtheorem{lemma}[theorem]{Lemma}
\newtheorem{corollary}[theorem]{Corollary}
\newtheorem{proposition}[theorem]{Proposition}

\theoremstyle{definition}
\newtheorem{definition}[theorem]{Definition}
\newtheorem{example}[theorem]{Example}

\newtheorem{problem}[theorem]{Problem}
\newtheorem{remark}[theorem]{Remark}

\newcommand{\C}{\mathbb{C}}
\newcommand{\D}{\mathbb{D}}
\newcommand{\N}{\mathbb{N}}

\renewcommand{\P}{\mathbb{P}}
\newcommand{\R}{\mathbb{R}}

\newcommand{\cA}{\mathcal{A}}

\newcommand{\cH}{\mathcal{H}}

\newcommand{\cO}{\mathcal{O}}
\newcommand{\cS}{\mathcal{S}}

\newcommand{\cT}{\mathcal{T}}

\newcommand{\cW}{\mathcal{W}}

\def\bs{\backslash}
\def\wt{\widetilde}

\def\a{\alpha}
\def\e{\epsilon}
\def\b{\beta}
\def\l{\lambda}

\numberwithin{equation}{section}

\begin{document}
\title 
[Holomorphic flexibility properties of complex manifolds]
{Holomorphic flexibility properties \\ of complex manifolds} 
\author{Franc Forstneri\v c}
\address{Institute of Mathematics, Physics and Mechanics, 
University of Ljubljana, Jadranska 19, 1000 Ljubljana, Slovenia}
\email{franc.forstneric@fmf.uni-lj.si}
\thanks{Research supported by grants P1-0291
and J1-6173, Republic of Slovenia.}

%
%
\subjclass [2000] {Primary 32E10, 32E30, 32H02, 32J25, 32Q28; secondary 14J99, 14R10}
\date{May 5, 2005} 
\keywords{Stein manifolds, holomorphic maps, algebraic manifolds, algebraic maps, 
Oka principle, transversality, dominability}

\begin{abstract}
We obtain results on approximation of holomorphic maps 
by algebraic maps, the jet transversality theorem for 
holomorphic and algebraic maps between certain classes of
manifolds, and the homotopy principle for holomorphic 
submersions of Stein manifolds to certain algebraic manifolds. 
\end{abstract}

\maketitle

\section{Introduction}
%
%
%
%

In the present paper we use the term {\em holomorphic flexibility property}
for any of several analytic properties of complex manifolds which are
opposite to Koba\-ya\-shi-Eisenman-Brody hyperbolicity,
the latter expressing {\em holomorphic rigidity}.
A connected complex manifold $Y$ is {\em $n$-hyperbolic} for some 
$n\in \{1,\ldots, \dim Y \}$ if every entire holomorphic 
map $\C^n\to Y$ has rank less than $n$; for $n=1$ this means that every 
holomorphic map $\C\to Y$ is constant  (\cite{Br}, \cite{Ei}, \cite{Ko1},
\cite{Ko2}). On the other hand, all flexibility properties of $Y$ 
will require the existence of many such maps.

We shall center our discussion around the following 
classical property which was studied by many authors 
(see the surveys \cite{Le} and \cite{F5}):

\smallskip
\textbf{Oka property}: 
{\em Every continuous map $f_0\colon X\to Y$ 
from a Stein manifold $X$ is homotopic 
to a holomorphic map $f_1\colon X\to Y$; if $f_0$ is holomorphic 
on (a neighborhood of) 
a compact $\cH(X)$-convex subset $K$ of $X$ then a homotopy 
$\{f_t\}_{t\in [0,1]}$ from $f_0$ to a holomorphic map $f_1$ 
can be chosen such that $f_t$ is holomorphic and 
uniformly close to $f_0$ on $K$ for every $t\in [0,1]$.} 
\smallskip

Here, $\cH(X)$-convexity means convexity with respect to the
algebra $\cH(X)$ of all holomorphic functions on $X$ (\S 2). 
It was recently proved in \cite{F6} that the 
Oka property of a complex manifold $Y$ is 
equivalent to the following

\smallskip
\textbf{Convex approximation property (CAP)}:
{\em Every holomorphic map $K\to Y$ on a compact convex set 
$K \subset \C^n$ $(n\in \N)$ can be approximated uniformly 
on $K$ by entire holomorphic maps $\C^n\to Y$}. 
\smallskip

In the present paper we consider relations
between the Oka property and various notions 
of {\em ellipticity} introduced by Gromov \cite{Gr2} and 
the author \cite{F2}, validity of the 
{\em jet transversality theorem} 
for holomorphic and algebraic maps 
(theorems \ref{t1.3}, \ref{t1.4} and results in \S 4),
and {\em dominability} by complex Euclidean spaces; 
these are summarized in corollary \ref{hierarchy}.
Another main result concerns approximation of holomorphic maps 
by algebraic maps (theorem \ref{t1.1}, corollary \ref{c1.2} 
and results in \S 3). 

A principal application is the 
{\em one-parametric homotopy principle for holomorphic submersions} 
from Stein manifolds to certain quasi-projective algebraic manifolds
(theorem \ref{t5.1}), extending the results from \cite{F3} and \cite{F4}.

An {\em algebraic manifold} (resp.\ an {\em algebraic variety}) 
will be understood in the sense of Serre (\S 34 in \cite{Se1}, p.\ 226): 
A space $X$ endowed with a Zariski topology 
(each decreasing sequence of closed sets is stationary),
together with a sheaf of rings $\cO_X$ of continuous $\C$-valued functions, 
such that $X$ is covered by finitely many open sets $U_j$,
each of them isomorphic (as a ringed space) to a quasi-affine 
algebraic manifold (resp.\ variety) 
$V_j\subset \C^{n_j}$ (Axiom VA$_I$), with the diagonal 
$\Delta$ of $X\times X$ closed in $X\times X$ (Axiom VA$_{II}$). 
Here, as usual, a {\em quasi-affine} (resp.\ {\em quasi-projective}) variety 
is the difference $X=X_0\bs X_1$ of two closed affine 
(resp.\ projective algebraic) subvarieties. 
Each algebraic manifold (resp.\ algebraic variety) $X$ 
has an underlying structure of a complex (holomorphic) manifold 
(resp.\ a complex space), with an induced map of $\cO_X$ 
to the sheaf $\cH_X$ of germs of holomorphic functions on $X$ 
\cite{Se2}. All algebraic maps in this paper are assumed 
to be morphisms (without singularities), thus defining holomorphic maps
of the underlying holomorphic manifolds.

The following notions are explained more precisely
in \S 2 below; see especially definition \ref{def1}.
A {\em spray} on a complex manifold $Y$ is a holomorphic 
map $s\colon E\to Y$ from the total space of a holomorphic 
vector bundle $p\colon E\to Y$, satisfying $s(0_y)=y$ 
for all $y\in Y$. The spray is algebraic if $p\colon E\to Y$ 
is an algebraic vector bundle and the spray map $s\colon E\to Y$ is algebraic. 
A complex (resp.\ algebraic) manifold $Y$ is (algebraically) 
{\em subelliptic} if it admits a finite collection of (algebraic) sprays 
$s_j\colon E_j\to Y$ such that for every $y\in Y$ the vector subspaces 
$(ds_j)_{0_y}(E_{j,y}) \subset T_y Y$ together span $T_y Y$; 
$Y$ is {\em elliptic} if this holds with a single 
(dominating) spray. Every complex homogeneous manifold
is elliptic \cite{Gr2} (see \S 2 below).

We begin by discussing algebraic approximations of holomorphic maps.

\begin{theorem} 
\label{t1.1}
If $X$ is an affine algebraic manifold and $Y$ is an algebraically 
subelliptic manifold then a holomorphic map $X\to Y$ which is homotopic 
to an algebraic map can be approximated by algebraic maps
uniformly on compacts in $X$. In particular, every null-homotopic 
holomorphic map $X\to Y$ is a limit of algebraic maps. 
\end{theorem}

Letting $K$ be a (geometrically) convex compact set
in $X=\C^n$ we obtain the following corollary.

\begin{corollary} 
\label{c1.2} 
{\bf (Algebraic CAP)}
If $Y$ is an algebraically subelliptic manifold then 
every holomorphic map $K \to Y$ from a compact convex subset 
$K\subset \C^n$ $(n\in \N)$ can be approximated uniformly on $K$ by 
algebraic maps $\C^n\to Y$. 
\end{corollary}

%
%
%
%
More precise results and examples can be found in \S 2 and \S 3.

The problem of approximating holomorphic maps by algebraic
maps is of central importance in analytic geometry. 
Algebraic approximations in general do not exist even 
for maps between affine algebraic manifolds (for example,
there are no nontrivial algebraic morphisms 
$\C\to \C_*=\C\bs\{0\}$).
Demailly, Lempert and Shiffman \cite{DLS} and Lempert \cite{Lm} proved 
that a holomorphic map from a Runge domain in an affine algebraic 
variety to a quasi-projective algebraic manifold can be approximated
uniformly on compacts by Nash algebraic maps. 
(A map $U \to Y$ from an open set $U$ in an 
algebraic variety $X$ is {\em Nash algebraic\/} if its graph 
is contained in an algebraic subvariety $\Gamma\subset X\times Y$ 
with $\dim \Gamma= \dim X$, \cite{Na}.) Nash algebraic 
approximations do not suffice in the proof of our theorem \ref{t5.1} 
where we need to approximate a holomorphic map 
$K \to Y$ on a compact convex set $K\subset \C^n$ by an entire 
map $\C^n\to Y$ whose ramification locus is a thin 
algebraic subvariety of $\C^n$. Global Nash algebraic maps 
would suffice for this purpose, but these are algebraic morphisms 
according to Serre (\cite{Se2}, p.\ 13, proposition~8).

%
%
%
%

We shall next discuss the {\em transversality theorems}
for holomorphic maps. 
If $X$ and $Y$ are smooth manifolds, $k\in \{0,1,2,\ldots\}$
and $Z$ is a smooth closed submanifold of $J^k(X,Y)$
(the manifold of $k$-jets of smooth maps of $X$ to $Y$)
then for a generic smooth map $f\colon X\to Y$ 
its $k$-jet extension $j^k f\colon X\to J^k(X,Y)$  
is transverse to $Z$ (Thom \cite{T1}, \cite{T2}; for extensions 
see \cite{Ab}, \cite{Go}, \cite{GM}, \cite{MT}, \cite{Tr}, \cite{Wh}). 
The analogous result only rarely holds for holomorphic maps. 
Indeed, the transversality theorem for one-jets 
of holomorphic maps $\C^n\to Y$ implies 
that all Kobayashi-Eisenman metrics on $Y$ vanish
identically, and $Y$ is {\em dominable} by $\C^n$, $n=\dim Y$. 
If such $Y$ is compact and connected then  its Kodaira dimension 
$\kappa= {\rm kod}\, Y$ (\cite{BH}, p.\ 29) satisfies 
$\kappa <  \dim Y$, i.e., $Y$ {\em is not of general Kodaira type} 
\cite{CG}, \cite{KO}, \cite{Kd1}, \cite{Kd2}. 

In the positive direction, the basic transversality 
theorem (for $0$-jets) holds for holomorphic maps to any manifold 
with a submersive family of holomorphic self-maps 
(Abraham \cite{Ab}); a classical example is {\em Bertini's theorem} to the 
effect that almost all projective hyperplanes in $\P_n= {\P_n}(\C)$ 
intersect a given complex submanifold $Z\subset \P_n$ 
transversally (\cite{GM}, p.\ 150; \cite{Kl}).
The jet transversality theorem holds for holomorphic maps of 
Stein manifolds to Euclidean spaces \cite{Fo}. 
Kaliman and Zaidenberg \cite{KZ} proved the jet transversality
theorem for holomorphic maps from Stein manifolds to {\em any} 
complex manifold provided that one shrinks the domain of the map 
(theorem \ref{t4.8} below).  

In \S 4 of this paper we prove the following transversality theorems.

%
%
%
%

%
%
\begin{theorem}
\label{t1.3}
Holomorphic maps from a Stein manifold to a subelliptic manifold 
satisfy the jet trans\-versality theorem with respect to closed
complex analytic subvarieties. Algebraic maps from an affine algebraic 
manifold to an algebraically subelliptic manifold satisfy the jet
transversality theorem on compact sets with respect 
to closed complex analytic (not necessarily algebraic) subvarieties. 
\end{theorem}

%
%
\begin{theorem}
\label{t1.4}
If $X$ is a Stein manifold and $Y$ is a complex manifold enjoying 
the Oka property then holomorphic maps $X\to Y$ satisfy 
the jet transversality theorem with respect to closed
complex analytic subvarieties.
\end{theorem}

These follow from theorems \ref{t4.2}, \ref{t4.3}
and proposition \ref{p4.6} in \S 4 where the reader can
find more precise formulations and further results.
Although the part of theorem \ref{t1.3} for holomorphic maps
is implied by theorem \ref{t1.4} and the fact that 
subellipticity implies the Oka property \cite{F2}, 
our  proof of theorem \ref{t1.3} is more elementary 
than the proof of the latter implication in \cite{F2}
and it also applies in the algebraic category where 
the Oka property is unavailable; for these reasons we 
separate them.

A property of a holomorphic map $X\to Y$ is said to be {\em generic} 
if it holds for all maps in a certain set of the second category in 
the Fr\`ech\'et space $\cH (X,Y)$ of all holomorphic maps $X\to Y$.
The jet transversality theorem implies that {\em singularities} of a 
generic  map $f\in \cH(X,Y)$ (points of nonmaximal rank) 
satisfy the codimension conditions in \cite{Fo}, proposition 2. 
This implies the following; for the last statement concerning 
injective immersions one needs
a multi-jet transversality theorem which is an easy extension 
(see \cite{Fo}, \S 1.3 for the case $Y=\C^N$).

\begin{corollary} 
\label{c1.5} 
If $X$ is a Stein manifold and $Y$ enjoys the Oka property then 
a generic holomorphic map $X\to Y$ is an immersion when 
$\dim Y\ge 2\, \dim X$ and an injective immersion when
$\dim Y\ge 2\, \dim X+1$. 
\end{corollary}

A holomorphic map $\pi \colon Y\to Y_0$ is called a 
{\em subelliptic Serre fibration} if it is a 
surjective subelliptic submersion (definition \ref{def1} below) 
and a Serre fibration, i.e., it enjoys the homotopy 
lifting property (\cite{W}, p.\ 8). Examples include holomorphic 
fiber bundles with subelliptic fibers and, more generally, 
subelliptic submersions which are topological fiber bundles 
(such as the unramified elliptic fibrations without exceptional
fibers; see \cite{BH}, p.\ 200). 
By theorem 1.8 in \cite{F6}, CAP {\em ascends 
and descends in a subelliptic Serre fibration}, 
in the sense that the manifolds $Y$ and $Y_0$ 
satisfy CAP at the same time. 
(See also \cite{Gr2}, 3.3.C' and 3.5.B'', and \cite{L2}, \cite{L3}.) 
For holomorphic fiber bundles the same conclusion
holds if the fiber satisfies CAP. A finite induction 
and the equivalence CAP$\Leftrightarrow$Oka property
imply the following result.

\begin{corollary} 
\label{c1.6}
Let $Y=Y_m\to Y_{m-1}\to \cdots\to Y_0$ where every map 
$Y_j\to Y_{j-1}$ $(j=1,2,\ldots,m)$ is a subelliptic  
Serre fibration. If one of the manifolds $Y_j$ enjoys the Oka property 
then all of them do. This holds in particular if $\dim Y_0=0$, 
and in such case $Y$ will be called {\em semisubelliptic}.  
\end{corollary}

If $\pi\colon Y\to Y_0$ is a {\em ramified holomorphic map}
then the Oka property of $Y_0$ need not imply the Oka property of $Y$. 
(A meromorphic function on a compact hyperbolic Riemann surface 
$Y$ is a finite ramified holomorphic map $Y\to Y_0=\P_1$, the Oka 
property holds for $\P_1$ but fails for $Y$. See also example 
\ref{Winkel1} in \S 6.) We don't know whether the Oka property 
of $Y$ implies the same for $Y_0$ (problem \ref{descending}).

%
%
%
%
Recall that a $p$-dimensional complex manifold $Y$ is 
(holomorphically) {\em dominable} by $\C^p$ if there exists a 
holomorphic map $f\colon \C^p\to Y$ with rank 
$p$ at $0\in \C^p$; if $Y$ and $f$ are algebraic then 
$Y$ is {\it algebraically dominable} by $\C^p$. 
If $Y$ is compact and connected then dominability 
by $\C^p$ implies $\mbox{kod}\,Y < p={\rm dim}\, Y$ \cite{Kd2}.

%
%

\begin{corollary}
\label{hierarchy}
{\bf (Hierarchy of holomorphic flexibility properties)}
The following implications hold for every complex manifold:
\begin{eqnarray*}
 \qquad {\rm homogeneous} \Longrightarrow 
 {\rm elliptic} \Longrightarrow {\rm subelliptic} \Longrightarrow 
 {\rm semisubelliptic}  \Longrightarrow \qquad \\ 
 {\rm CAP}\Longleftrightarrow   {\rm Oka\ property} 
  \Longrightarrow {\rm jet\  transversality\ theorem} 
  \Longrightarrow {\rm dominable}. 
\end{eqnarray*}
In the algebraic category we have the implications 
$$
\begin{array}[t]{*5c}   
            {\rm elliptic} & \Longrightarrow &  {\rm subelliptic} 
            & \Longrightarrow & {\rm CAP}   \\
            & & \Downarrow & & \Downarrow \\
             & & {\rm jet\ transversality} 
             & \Longrightarrow & {\rm dominable}. \\
\end{array}
$$
\end{corollary}

The `jet transversality theorem' refers to holomorphic maps
from any Stein manifold to $Y$. In the algebraic category 
{\rm CAP} is interpreted in the sense of corollary \ref{c1.2}, 
and the first vertical arrow in the last display is theorem \ref{t1.3}. 

Clearly CAP of $Y$ implies the following {\em strong dominability}: 
For every point $y\in Y$ there is a holomorphic 
map $f_y\colon \C^p\to Y$ $(p=\dim Y)$ such that $f_y(0)=y$ 
and $df_y$ has maximal rank $p$ at $0\in\C^p$. 
(Note that a family $\{f_y\}_{y\in Y}$ of such maps which depend  
holomorphically on the base point $y\in Y$ is a dominating spray on $Y$.) 
If $Y$ is strongly dominable by $\C^p$ $(p=\dim Y)$ then every 
bounded plurisubharmonic function on $Y$ is constant.  
If such $Y$ is also Stein then every $\R$-complete holomorphic 
vector field on $Y$ is also $\C$-complete, i.e., it induces a holomorphic 
action of $(\C,+)$ (corollary 2.2 in \cite{F1}). 
This implies the following.

\begin{corollary} 
\label{completefields}
If a Stein manifold $Y$ enjoys {\rm CAP} then every 
$\R$-complete holomorphic vector field on $Y$ is also $\C$-complete.
\end{corollary}

%
%
\begin{remark} 
\label{generalOka} 
{\bf (CAP and Oka-type properties)}
A more precise term for the Oka property as defined above is
{\em the basic Oka property with approximation}, where
`basic' refers to the non-parametric case.
More general Oka-type properties have been studied
by many authors 
(see e.g.\ \cite{G1}, \cite{G2}, \cite{G3}, \cite{Ca}, \cite{GK}, 
\cite{FR}, \cite{Le}, \cite{Gr2}, \cite{HK}, \cite{Wi},
\cite{HL}, \cite{FP1}, \cite{FP2}, \cite{FP3}, \cite{F2},
\cite{L2}, \cite{L3}). 
The {\em basic Oka property with (jet) interpolation}
of $Y$ refers to the possibility of homotopically 
deforming a continuous map $f_0\colon X\to Y$ from any Stein  
manifold $X$ to a holomorphic map $f_1\colon X\to Y$, 
keeping it fixed (to any finite order in the jet case) 
on a closed complex subvariety $X_0\subset X$ along 
which $f_0$ is holomorphic. Combining interpolation
with approximation on compact $\cH(X)$-convex subsets one
obtains the {\em Oka property with (jet) interpolation and approximation}.
By corollary 1.4 in \cite{F7} all these properties are equivalent 
to CAP, and hence 
{\em the conclusion of corollary \ref{c1.6} above holds for 
all mentioned Oka-type properties.} The analogous equivalences 
hold for the {\em parametric Oka properties} 
(theorem 5.1 in \cite{F6} and   theorem 6.1 in \cite{F7}). 
All these Oka properties (also the parametric versions) 
are implied by ellipticity \cite{Gr2} and subellipticity \cite{F2}, 
and they are equivalent to ellipticity on a Stein manifold 
(\cite{Gr2}, 3.2.A.; proposition 1.2 in \cite{FP3}).
\end{remark}

\begin{remark}
\label{LP}
As pointed out in \cite{F6}, CAP is a localization 
($=$restriction) of the Oka property 
to model pairs --- $K$ a compact convex set in $X=\C^n$. 
The equivalence CAP$\Leftrightarrow\,$Oka property
supports the following heuristic principle; 
compare with the formulation of the {\em Oka principle} 
by Grauert and Remmert (\cite{GRe}, p.\ 145).

\smallskip
\textbf{Localization principle}: 
{\em
Problems concerning holomorphic maps from Stein manifolds 
have only homotopical obstructions provided that the target 
manifold satisfies a suitable flexibility property expressed
in terms of holomorphic maps from Euclidean spaces.}
\smallskip

This principle, once established for a certain class of maps,
can substantially simplify the analysis in concrete examples
as is amply demonstrated in  \cite{F6} and in \S 6 below. 
On the philosophical level it reveals the interesting but
not entirely surprisig fact that many complex analytic problems 
involving Stein manifolds reduce to problems on Euclidean
spaces (which does not necessarily make them easy).

By \cite{F6} the localization principle holds for sections 
of holomorphic fiber bundles over any Stein manifold, 
the corresponding local property being CAP of the fiber. 
Many classical problems on Stein manifolds, including the 
classification of principal holomorphic fiber bundles 
and their associated bundles with homogeneous fibers,
fits into this framework (see Grauert \cite{G3}, Cartan \cite{Ca}, 
and the survey \cite{Le} of Leiterer).

One may also consider special classes of holomorphic maps such 
as immersions, submersions, maps of constant rank, etc.
According to \cite{F4} {\em the localization principle holds for 
holomorphic submersions} $X\to Y$ from Stein manifolds $X$,
the corresponding local property of $Y$ being Property 
${\mathrm S}_{\mathrm n}$ with $n=\dim X \ge \dim Y$
(see definition \ref{def3} and theorem \ref{t5.1} in \S 5 below).
For {\em holomorphic immersions} a localization principle 
is not known at this time, but the {\em homotopy principle 
for immersions of Stein manifolds to Euclidean spaces} was proved 
by Eliashberg and Gromov \cite{EG0}, \cite{Gr1}.
\end{remark}

The flexibility properties discussed in this paper are mutually
exclusive with any kind of hyperbolicity property 
(for the latter see e.g.\ \cite{Br}, \cite{Ei}, \cite{Ko1}, \cite{Ko2}, 
\cite{Za} and the surveys \cite{D2}, \cite{D3}, \cite{PS}, \cite{Si2}. 
In particular, none of the flexibility properties holds 
for compact complex manifolds of Kodaira general type;
it is conjectured that every such manifold $Y$ is almost
hyperbolic in the sense that there is a proper 
complex subvariety $Y_0\subset Y$ containing the image of any nonconstant 
holomorphic map $\C\to Y$ (the Green-Griffiths conjecture, \cite{GG}).
In \S 6 we consider flexibility properties of complex curves 
and surfaces; in the latter case our analysis rests upon the 
fairly complete results of Buzzard and Lu \cite{BL} on 
dominability of compact complex surfaces with 
less than maximal Kodaira dimension.

%
%
%
\section{Preliminaries}
We denote by $\cH _X$ the sheaf of germs of holomorphic functions 
on a complex (holomorphic) manifold $X$, and by $\cH (X)=\cH (X,\C)$ the 
algebra of all global holomorphic functions on $X$. Given a pair of complex 
manifolds $X,Y$, we denote by $\cH (X,Y)$ the set of all holomorphic maps 
$X\to Y$, endowed with the compact-open topology. 
This topology is induced by a complete metric and hence 
$\cH (X,Y)$ is a {\em Baire space}. 
A property of $f\in \cH (X,Y)$ is said to be 
{\em generic} if it holds for all $f$ in a {\em residual subset\/} 
(a countable intersection of open dense subsets) 
of $\cH(X,Y)$.

A function (or a map) is said to be {\em holomorphic on 
a compact set} $K$ in a complex manifold $X$ if 
it is holomorphic in an open set $U\subset X$ containing $K$.
A homotopy of maps $\{f_t\}$ is holomorphic on  $K$ if there is
an open  neighborhood $U$ of $K$, independent of $t$, 
such that $f_t$ is holomorphic on $U$ for every $t$. 

An {\it affine manifold\/} is a closed complex or algebraic 
submanifold of a complex Euclidean space; an affine 
complex manifold is the same thing as a {\em Stein manifold}
according to the embedding theorems \cite{Bi}, \cite{EG},
\cite{Ns},  \cite{Re}, \cite{Sch}. 

A compact subset $K$ in a Stein manifold $X$ 
is $\cH (X)$-convex ({\it holomorphically convex} in $X$) if 
for any $p\in X\backslash K$ there exists $f\in \cH (X)$ with 
$|f(p)|>\sup_K |f|$. If $X$ is embedded in $\C^N$ then a set $K\subset X$ 
is $\cH (X)$-convex if and only if it is 
$\cH (\C^N)$-convex, i.e., {\it polynomially convex}.

We denote by $\cO_X$ the structure sheaf of an algebraic 
manifold $X$ and by $\cO (X)$ the algebra of all regular algebraic 
functions on $X$. As we have already said in \S 1,
an algebraic manifold or variety will be understood in the sense
of Serre (\S 34 in \cite{Se1}, p.\ 226); in particular,
there is an underlying complex manifold (resp.\ analytic variety) 
structure on $X$. By $\cO(X,Y)$ we denote the set of all 
regular algebraic maps (morphisms) between a pair of algebraic manifolds.
Clearly $\cO(X,Y) \subset \cH (X,Y)$ but $\cO (X,Y)$ need not be 
closed in $\cH (X,Y)$. If both $X$ and $Y$ are projective
algebraic then $\cO(X,Y)=\cH(X,Y)$ by Sere's {\em GAGA principle} 
\cite{Se2}.

%
%
%
%
A {\em fiber-spray} associated to a holomorphic 
submersion $h\colon Y\to Y'$ is a triple $(E,p,s)$ consisting of a 
holomorphic vector bundle $p\colon E\to Y$ and a holomorphic map 
$s\colon E\to Y$ such that for each $y\in Y$ we have $s(0_y)=y$ 
and $s(E_y) \subset Y_{h(y)}= h^{-1}(h(y))$
(see \cite{Gr2}, \S 1.1.B., and \cite{FP1}).
A spray on a complex manifold $Y$ is a fiber-spray associated 
to the trivial submersion $Y\to point$. For each $y\in Y$ let 
$VT_y Y =\ker dh_z \subset T_y Y$, the 
{\em vertical tangent space} of $Y$ with respect to $h$. 

%
%
\begin{definition}
\label{def1}
A holomorphic submersion $h\colon Y\to Y'$
is {\em subelliptic} if each point in $Y'$ has an open neighborhood 
$U\subset Y'$ such that the restricted submersion
$h\colon Y|_U=h^{-1}(U) \to U$ admits finitely 
many fiber-sprays $(E_j,p_j,s_j)$ $(j=1,\ldots,k)$ satisfying  
$$ 
	(ds_1)_{0_y}(E_{1,y}) + (ds_2)_{0_y}(E_{2,y})\cdots 
                     + (ds_k)_{0_y}(E_{k,y})= VT_y Y   \eqno(2.1)
$$
for each $y\in Y|_U$; such a collection of sprays 
is said to be {\it dominating}. The submersion is {\em elliptic} 
if the above holds with $k=1$. A complex manifold $Y$ is 
(sub-)elliptic if the trivial submersion $Y\to point$ is such. 
An algebraic manifold $Y$ is {\em algebraically subelliptic} if 
it admits finitely many algebraic sprays satisfying the 
domination condition (2.1) for the trivial submersion 
$Y\to point$; $Y$ is {\em algebraically elliptic} if 
it admits a dominating algebraic spray. 
\end{definition}

Examples of (sub-)elliptic manifolds and submersions can 
be found in \cite{G3} (especially sections 0.5.B and 3.4.F), 
\cite{F2}, \cite{FP1}. If a complex Lie group $G$, 
with Lie algebra ${\mathbf g}$,
acts holomorphically and transitively on a 
complex manifold $Y$ then the map $Y\times {\mathbf g}\to Y$,
$(y,v)\to e^v y$, is a dominating spray on $Y$. 
Furthermore, {\em every algebraic Lie group $G$ which admits
no homomorphisms to $\C_*$ is algebraically elliptic};
the following proof has been explained to me by J.\ Winkelmann.
The condition on $G$ implies that it is generated by its unipotent subgroups.
Each left invariant vector field on $G$ arising from a
unipotent subgroup gives rise to an algebraic spray on $G$.
Under the above generation condition these vector fields 
span the tangent space of $G$ at each point, and the composition 
of their flows gives a dominating algebraic spray on $G$. 
The Lie group $\C_*$ is not algebraically elliptic.

The following result from 
\cite{F2} (which is implicitly contained in 
lemmas 3.5.B.\ and 3.5.C.\ in \cite{Gr2})
will be important for us.

%
%
\begin{proposition}
\label{p2.2}
Let $Y$ be a quasi-projective algebraic manifold. 
If every point of $Y$ has a Zariski open neighborhood
$U\subset Y$ such that $U$ is algebraically subelliptic
then $Y$ is itself algebraically subelliptic.
\end{proposition}

The following classes of algebraic manifolds will be used in \S 5.

%
%
\begin{definition}
\label{def2}
Let $Y$ be a quasi-projective algebraic manifold.
\begin{itemize}
\item[(a)] $Y$ is of {\em Class} $\cA_0$ 
if it is covered by finitely many Zariski open sets
biregularly isomorphic to the affine space $\C^p$
with $p=\dim Y$. 
\item[(b)] $Y$ is of {\em Class} $\cA$ if 
$Y=\widehat Y\bs A$ where $\widehat Y$ is of class $\cA_0$ 
and $A$ is a thin ($=$of codimension at least two)
algebraic subvariety of $\widehat Y$. 
\end{itemize}
\end{definition}

Thus a $p$-dimensional manifold of Class $\cA$ is Zariski locally 
equivalent to $\C^p\bs A$ where $A$ is a thin 
algebraic subvariety in $\C^p$, depending on the chosen
Zariski open set in $Y$. Such manifold $\C^p\bs A$ is algebraically elliptic 
(\cite{FP1}, p.\ 119; \cite{Gr2}). Indeed, there is a 
polynomial spray $s\colon \C^p\times\C^N \to \C^p$ 
which maps $(\C^p\bs A) \times\C^N$ to $\C^p\bs A$ 
and is dominating at every point $\C^p\bs A$;
such $s$ is obtained by composing flows of suitably chosen 
algebraic shear vector fields on $\C^p$ which 
vanish on $A$ and generate the tangent bundle at each point
of $\C^p\bs A$. Hence proposition \ref{p2.2} implies

\begin{corollary}
\label{c2.4}
Every manifold of Class $\cA$ is algebraically subelliptic.
\end{corollary}

Class $\cA_0$ includes all complex affine and projective 
spaces, as well as Grassmanians. Another example 
is the total space $W$ of a holomorphic fiber bundle 
$\pi \colon W\to Y$ where the base $Y$ is a manifold 
of Class $\cA _0$, the fiber is $F=\pi^{-1}(y) \in \{\C^m,\P_m\}$, 
and the structure group is $GL_m(\C)$ respectively $PGL_m(\C)$. 
Every such bundle is algebraic by GAGA \cite{Se2}, 
and its restriction to any affine Zariski open set 
$\C^p \simeq U \subset Y$ is algebraically trivial, 
$\pi^{-1}(U)\simeq U\times F$. Hence $W$ is 
covered by Zariski open sets biregularly isomorphic to 
$\C^{p+m}$, i.e., $W$ is of Class $\cA_0$.  
An example are the {\em Hirzebruch surfaces} 
$H_l$, $l=0,1,2,\ldots$ (\cite{BH}, p.\ 191); 
these are $\P_1$-bundles over $\P_1$, and each of  
them is birationally equivalent to $\P_2$.

Manifolds of Class $\cA$ have been considered by 
Gromov who called them  {\it Ell-regular\/} and showed 
that this class is stable under blowing up points 
(\cite{Gr2}, \S 3.5.D''):

\begin{proposition} 
\label{p2.5}  
{\rm (\cite{Gr2})}  
If $Y$ is of Class $\cA$ (respectively of Class $\cA_0$) then
any manifold obtained from $Y$ by blowing up finitely 
many points is also of Class $\cA$ (respectively of Class $\cA _0$).
\end{proposition}

\begin{proof}  By localization it suffices to show that
the manifold $L$, obtained by blowing up $\C^q$ at the origin, 
is of Class $\cA $. $L$ is the total space of a holomorphic 
line bundle $\pi \colon L\to\P_{q-1}$ (the universal bundle); 
$L$ is trivial over the complement of each hyperplane 
$\P_{q-2}\subset \P_{q-1}$ (which equals $\C^{q-1}$), 
and hence every point in $L$ has a Zariski neighborhood 
of the form $\pi^{-1}(\P_{q-1}\backslash \P_{q-2})$ 
which is biregular to $\C^q$. 
\end{proof}

%
%
%
%
\section{Algebraic approximation}
All algebraic maps are assumed to be morphisms (without singularities).

%
%
\begin{theorem}
\label{t3.1} 
Let $X$ be an affine algebraic manifold and $K$ a compact 
$\cH (X)$-convex subset of $X$. Let $Y$ be an algebraically 
subelliptic manifold (def.\ \ref{def1}) with a distance function
$d$ induced by a Riemannian metric. If $f_t\colon K\to Y$ 
$(t\in [0,1])$ is a homotopy of holomorphic maps 
such that $f_0$ extends to an algebraic map $X\to Y$ then 
for every $\e>0$ there exists an algebraic map 
$F\colon X\times \C\to Y$ such that $F(\cdotp,0)=f_0$ and 
$d(F(x,t),f_t(x))<\e$ for every $x\in K$ and $t\in [0,1]$.
\end{theorem}

%
%
\begin{corollary}
\label{c3.2} 
If $K\subset X$ and $Y$ are as in theorem \ref{t3.1} then
every null-homotopic holomorphic map $K\to Y$ can be approximated 
uniformly on $K$ by algebraic maps $X\to Y$. In particular, 
if $K$ is a compact convex set in $\C^n$ then 
every holomorphic map $K\to Y$ can be approximated uniformly on 
$K$ by algebraic maps $\C^n\to Y$.
\end{corollary}

Examples of algebraically subelliptic manifolds are given by 
propositions \ref{p2.2}, \ref{p2.5} and corollary \ref{c2.4}. 
We do not know whether every homotopy class of maps $X\to Y$ 
in theorem \ref{t3.1} contains an algebraic map.

Theorem \ref{t3.1} is a special case (with $Z=X\times Y\to X$)
of the following.

\begin{theorem}
\label{t3.3} 
Let $h\colon Z\to X$ be an algebraic submersion 
from an algebraic manifold $Z$ onto an 
affine algebraic manifold $X$. Let $d$ be a distance function
on $Z$ induced by a Riemannian metric. Assume that $Z$ 
admits a family of algebraic fiber-sprays $(E_j,p_j,s_j)$ 
$(j=1,\ldots,k)$ satisfying the domination property (2.1)
at every point $z\in Z$.
Let $K\subset X$ be a compact $\cH (X)$-convex set
and let $f_t\colon K\to Z$ $(t\in [0,1])$ be a homotopy of 
holomorphic sections such that $f_0$ extends to an
algebraic section $X\to Z$. For every $\e>0$ 
there is an algebraic map $F\colon X\times \C\to Z$ such that
$h(F(x,t))=x$ for $(x,t)\in X\times\C$,
$F(\cdotp, 0)=f_0$, and $d(F(x,t),f_t(x))<\e$ 
for every $(x,t) \in K\times [0,1]$.
\end{theorem}

An algebraic submersion satisfying the hypotheses of 
theorem \ref{t3.3} will be called {\em algebraically subelliptic} 
(compare with def.\ \ref{def1}). Theorem \ref{t3.3} is an 
algebraic analogue of theorem 3.1 in \cite{F2}; see also \cite{Gr2} 
and theorems 4.1 and 4.2 in \cite{FP1}.

\smallskip
\noindent 
\textit{Proof of theorem \ref{t3.3}.}
Let $(E,p,s)$ be the  composed (algebraic) fiber-spray on $Z$
obtained from the fiber-sprays $(E_j,p_j,s_j)$
$(j=1,\ldots,k)$ (see \cite{Gr2}, \S 1.3, and \cite{FP1}, definition 3.3).
We briefly recall the construction. 
Set $(E^{(1)},p^{(1)},s^{(1)})=(E_1,p_1,s_1)$.
If $k>1$, let $E^{(2)} = s^{(1)*} (E_2) \to E^{(1)}$ 
(the pull-back of $p_2\colon E_2\to Z$ by the spray map 
$s^{(1)} \colon E^{(1)} \to Z$), and define 
$p^{(2)}\colon E^{(2)} \to Z$ and  
$s^{(2)} \colon E^{(2)} \to Z$ by 
$p^{(2)}(e_1,e_2)=p_1(e_1)$, $s^{(2)}(e_1,e_2)=s_2(e_2)$. 
(Note that $s_1(e_1)=p_2(e_2)$.) Continuing inductively 
we obtain a sequence of algebraic vector bundle projections
$$ 
	E=E^{(k)} \,\longrightarrow{}\, 
	\cdots \,\longrightarrow{}\, E^{(2)} \,\longrightarrow{}\,  
	E^{(1)} \,\longrightarrow{}\, Z. 
$$
The composed bundle $E\to Z$ with fibers $E_z$ ($z\in Z$)
has a well defined zero section which we identify with $Z$, and 
$TE|_Z\simeq TZ\oplus E_1\oplus\cdots\oplus E_k$.
Denote by $s=s^{(k)}\colon E\to Z$ the composed spray.
The restriction of its differential 
$ds_{0_z} \colon T_{0_z} E\to T_z Z$ to the fiber $E_{j,z}$
of $E_j$ over $z$ in the above direct sum decomposition equals 
$(ds_j)_{0_z} \colon T_{0_z} E_{j,z} \to VT_z Z$.
Hence (2.1) is equivalent to the domination 
property of the composed spray:  
$$ 
	(ds)_{0_z}(T_{0_z} E_z)=VT_z Z, \quad z\in Z. \eqno(3.1)
$$
The bundle $E\to Z$ admits a (noncanonical) holomorphic 
vector  bundle structure over any Stein subset of $Z$ 
(\cite{Gr2}, \S 1.3, and \cite{FP1}, corollary 3.5).

In the remainder of the proof we shall only work with 
the composed spray bundle $(E,p,s)$ and will no longer
need the individual sprays $(E_j,p_j,s_j)$.

\begin{lemma} 
\label{l3.4} 
Let $V\subset\subset U$ be open Stein neighborhoods of $K$ in $X$ 
such that the homotopy $f_t$ $(t\in [0,1])$ in theorem \ref{t3.3} 
is defined in $U$. Set $V_t=f_t(V) \subset Z$ for $t\in [0,1]$. 
There are numbers $l\in \N$ and $0=t_0<t_1< \cdots <t_l =1$
such that for every $j=0,1,\ldots,l-1$ there exists a homotopy of 
holomorphic sections $\xi_t$ of the restricted bundle 
$E|_{V_{t_j}} \to V_{t_j}$  $(t\in I_j=[t_j,t_{j+1}])$ such that 
$\xi_{t_j}$ is the zero section and $s(\xi_t(z))= f_t(h(z))$
for all $t\in I_j$ and $z\in V_{t_j}$.
\end{lemma}

\begin{proof}
Assume first that there exists a Stein open set
$\Omega\subset Z$ containing $\cup_{t\in[0,1]}\overline V_t$.
By corollary 3.5 in \cite{FP1} the restriction $E|_\Omega\to\Omega$ 
admits the structure of a holomorphic vector bundle and 
a holomorphic direct sum splitting $E|_\Omega=H\oplus H'$,
where $H'$ is the kernel of $ds$ at the zero section of $E$
(\cite{GRo}, p.\ 256, theorem 7). It follows from (3.1) 
that for every $z\in \Omega$ the 
restriction $s\colon H_z\to Z_{h(z)}=h^{-1}(h(z))$ 
maps a neighborhood of $0_z$ in $H_z$ biholomorphically onto a 
relative neighborhood of $z$ in the fiber $Z_{h(z)}$. 
The size of this neighborhood, and of its image in the 
corresponding fiber of $Z$, can be chosen uniform for points 
in the compact set $\cup_{t\in [0,1]} \overline V_t \subset\Omega$. 
Hence there is a $\delta >0$ such that for every $t\in [0,1]$ 
the local inverse of $s\colon H|_{V_t} \to Z$ 
at the zero section gives a homotopy of sections 
$\xi_\tau$ of $H|_{V_t}$ ($\tau\in J_t=[t,t+\delta] \cap [0,1]$),
with $\xi_t$ the zero section, such that $s(\xi_\tau(z))=f_\tau(h(z))$
for $\tau\in J_t$ and $z\in V_t$. This proves lemma \ref{l3.4}
in the special case. 

For the general case observe that $f_t(U)$, being a 
closed Stein submanifold of $Z|_U$, admits an open Stein 
neighborhood in $Z$ \cite{Si1}, \cite{D1}. By compactness there are 
Stein open sets $\Omega_j\subset Z$ $(j=1,2,\ldots,m)$ and a 
partitition $[0,1]= \cup_{j=1}^m I_j$ into adjacent closed 
subintervals $I_j$ such that $\cup_{t\in I_j} \overline V_t \subset \Omega_j$.
It remains to apply the above argument separately for each $I_j$. 
\end{proof}

Denote by $(E^{(l)},p^{(l)},s^{(l)})$ the $l$-th iterated
bundle of $(E,p,s)$ (definition 3.3 in \cite{FP1}, p.\ 132;
this is just the $l$-tuple composition of $(E,p,s)$ with itself).
Let $(E',p',s')$ denote the pull-back of $(E^{(l)},p^{(l)},s^{(l)})$ 
to $X$ by the algebraic map $f_0\colon X\to Z$; by the construction
this is an algebraic composed spray bundle over $X$.

\begin{lemma} 
\label{l3.5} 
There is a homotopy $\eta_t\colon V \to E'|_V$ $(t\in [0,1])$ of 
holomorphic sections of the restricted bundle 
$E'|_V\to V$ such that $\eta_0$ is the zero 
section and $s'(\eta_t(x))=f_t(x)$ for every $x\in V$ and $t\in [0,1]$.
\end{lemma}

\begin{proof} 
It suffices to assemble the individual homotopies 
$\{\xi_t\colon t \in [t_j,t_{j+1}]\}$ $(j=0,\ldots,l-1)$, 
furnished by lemma \ref{l3.4},
into a homotopy of sections $\wt \xi_t \colon V_0\to E^{(l)}|_{V_0}$ 
$(t\in [0,1])$ of the iterated bundle $E^{(l)}$ over the open subset 
$V_0= f_0(V)$ in the algebraic submanifold $f_0(X)$ of $Z$. 
Clearly $\wt \xi_t$ corresponds to a holomorphic sections 
$\eta_t\colon V\to E'=f_0^* (E^{(l)})$ of the pull-back bundle,
with $\eta_0$ being the zero section. 
(For further details see \cite{FP1}, 
proposition 3.6 on p.\ 134.)
\end{proof}

\begin{lemma}
\label{l3.6}
Let $d'$ be a distance function on $E'$ induced by a 
Riemannian metric. Let $\{\eta_t\}_{t\in [0,1]}$ 
be as in lemma \ref{l3.5}.
For every $\delta >0$ there is an algebraic map 
$\eta' \colon X\times \C \to E'$ satisfying 
\begin{itemize}
\item[(i)]  $\eta'(x,0)=0_x\in E'_x$ $(x\in X)$, and
\item[(ii)] $d'(\eta'(x,t),\eta_t(x)) < \delta$ for all
$x\in K$ and $t\in [0,1]$. 
\end{itemize}
\end{lemma}

\begin{proof} 
By construction of the composed bundle $E'\to X$ 
there is a finite sequence
$$
	E' = E^{m,0} \,\longrightarrow{}\, E^{m-1,0} \,\longrightarrow{}\, \cdots 
	\,\longrightarrow{}\, E^{1,0} \,\longrightarrow{}\,  X,  
								\eqno(3.2)
$$
with $m=kl$ and $E^{(1,0)}=f_0^* E_1 \to X$, in which every map 
$E^{j,0}\to E^{j-1,0}$ is an algebraic vector bundle projection. 
Here $k$ is the number of the initial sprays in theorem \ref{t3.3},
and $l$ is the number in lemma \ref{l3.4}.

Since $X$ is an affine algebraic manifold, 
the algebraic vector bundle $E^{1,0}\to X$ is generated by 
finitely many (say $n_1$) algebraic sections
according to Serre's theorem A (\cite{Se1}, p.\ 237, Th\'eor\`eme 2).
This gives a surjective algebraic map 
$\pi_1 \colon E^{1,1} = X\times \C^{n_1}\to E^{1,0}$
of the trivial rank $n_1$ bundle onto $E^{1,0}$. Pulling back
the sequence (3.2) to the new total space $E^{1,1}$ 
we obtain a commutative diagram
$$  
\begin{array}[t]{*{11}c}
        E^{m,1} & \longrightarrow{} & E^{m-1,1} & \longrightarrow{}  &
         \cdots & \longrightarrow{} & E^{2,1} & 
         \longrightarrow{} & E^{1,1} &  \longrightarrow{} & X \\
	  \;\;\downarrow{\pi_m} && \;\;\downarrow{\pi_{m-1}} &&&& \;\;\downarrow{\pi_2} && 
	  \;\;\downarrow{\pi_1} & & \Vert \\
	 E^{m,0} & \longrightarrow{} & E^{m-1,0} & \longrightarrow{} & 
	 \cdots & \longrightarrow{} & E^{2,0} & 
	 \longrightarrow{} & E^{1,0} & \longrightarrow{} & X  
\end{array}
$$
in which all horizontal maps are algebraic vector bundle projections
and the vertical maps $\pi_j$ for $j\ge 2$ are the induced natural maps
which are bijective on fibers.
More precisely, we begin by letting $E^{2,1} \to E^{1,1}$ be the pull-back
of the vector bundle $E^{2,0} \to E^{1,0}$ (in the bottom row) by the 
vertical morphism $\pi_1\colon E^{1,1}\to  E^{1,0}$, and we denote
by $\pi_2 \colon E^{2,1} \to E^{2,0}$ the asociated natural map
which makes the respective diagram commute. Moving one step to 
the left, $E^{3,1} \to E^{2,1}$ is the pull-back
of the bundle $E^{3,0} \to E^{2,0}$ in the bottom
row by the vertical morphism $\pi_2 \colon E^{2,1}\to  E^{2,0}$,
and $\pi_3\colon E^{3,1} \to E^{3,0}$ is the associated natural map; etc.
There is an algebraic spray map $s^1 \colon E^{m,1} \to Z$ which 
is the composition of $\pi_m\colon E^{m,1}\to E^{m,0}$ 
with the initial spray $s\colon E^{m,0}=E \to Z$.

We claim that the homotopy $\eta_t$ of holomorphic sections
of $E^{m,0}|_V=E'|_V \to V$, furnished by lemma \ref{l3.5},
lifts to a homotopy  $\eta_t^1$ of holomorphic 
sections of $E^{m,1}|_V \to V$ such that $s^1(\eta_t^1)=f_t$ for 
all $t\in [0,1]$, and $\eta_0^1$ is the zero section.
It suffices to see that the $E^{1,0}$-component 
of $\eta_t$ (i.e., the projection of $\eta_t$ under the composed 
projection $E^{m,0}\to E^{1,0}$) lifts to $E^{1,1}$; the rest of the
lifting is then obtained by applying the inverses 
of the fiberwise isomorphic vertical maps.
But this follows from the fact that the surjective 
vector bundle map $\pi_1\colon E^{1,1}\to E^{1,0}$
admits a holomorphic splitting $\sigma_1\colon E^{1,0} \to E^{1,1}$ 
over $X$, with $\pi_1\circ\sigma_1$ the identity on $E^{1,0}$
(theorem 7 in \cite{GRo}, p.\ 256).

Repeating the same argument with the bundle 
$E^{2,1}\to E^{1,1}=X_1$ over the affine manifold
$X_1=X\times \C^{n_1}$ we obtain a surjective algebraic vector bundle map 
$E^{2,2}=X_1\times \C^{n_2}=X\times \C^{n_1+n_2} \to E^{2,1}$. 
As before we lift the top line in the above diagram to a new level 
$$
	E^{m,2} \,\longrightarrow{}\, E^{m-1,2} \,\longrightarrow{}\, 
	\cdots \,\longrightarrow{}\, E^{2,2} \,\longrightarrow{}\, 
	E^{1,1} = X_1=X\times\C^{n_1}.
$$
Note that $E^{2,2}=X_1\times\C^{n_2} = X\times \C^{n_1+n_2}$ 
(algebraic equivalence). The homotopy $\eta_t^1$ lifts to a homotopy
$\eta^2_t \colon V\to E^{m,2}|_V$, with $\eta^2_0$ the zero section,
and we have a new spray map $s^2\colon E^{m,2}\to Z$ 
satisfying $s^2(\eta^2_t)=f_t$ for all $t\in [0,1]$. 

Continuing inductively we obtain after $m$ steps 
a lifting of the homotopy $\eta_t$ to a homotopy 
$\eta_t^m \colon V\to E^{m,m}|_V =V\times \C^N$ 
($t\in [0,1]$), consisting of holomorphic sections of 
$E^{m,m}=X\times \C^N$ ($N=n_1+n_2+\cdots +n_m$) over 
the open subset $V \subset X$, with $\eta_0^m$ being the zero section.
By construction there is an algebraic spray 
$s^m\colon E^{m,m}\to Z$ such that 
$s^m(\eta_t^m)=f_t \colon  V\to Z|_V$ for all $t\in [0,1]$. 

Recall that $X$ is a closed algebraic submanifold of an affine
space $\C^n$. The $\cH(X)$-convex set $K \subset X$ is then 
polynomially convex in $\C^n$, and $K\times [0,1]$ is polynomially convex
in $\C^{n+1}$. (We have identified the segment $[0,1] \subset\R$
with its image in $\C$.) By a small extension of 
the Oka-Weil theorem we can approximate
the homotopy $\{\eta_t^m\}_{t\in[0,1]}$ 
(which is continuous in $(x,t) \in V\times [0,1]$ and 
holomorphic with respect to $x \in V$ 
for every fixed $t \in [0,1]$) uniformly on 
the set $K\times [0,1]$ by a holomorphic polynomial map 
$\C^n\times\C\to\C^n\times\C^N$ of the form 
$\wt g(x,t)= (x, g(x,t))$, with
$g(x,0)=0$ for $x\in \C^n$. By projecting the 
point $\wt g(x,t) \in E^{m,m}=X\times \C^N$ 
$(x\in X,\ t\in \C)$ back to 
the bundle $E'=E^{m,0}$ we obtain an algebraic map 
$\eta'(x,t)$ satisfying lemma \ref{l3.6}. 
\end{proof}

If $s'$, $\eta'$ and $\delta$ are as in lemmas 
\ref{l3.5} and \ref{l3.6}, with $\delta>0$ chosen 
sufficiently small, then the algebraic map 
$$
	F(x,t)=  s'(\eta'(x,t)) \in Z, \quad 
	(x,t)\in X\times \C
$$
satisfies the conclusion of theorem \ref{t3.3}.

%
%
%
%
\section{Transversality theorems for holomorphic and algebraic maps}
Given a pair of complex manifolds $X$, $Y$ we denote 
by $J^k(X,Y)$ the manifold of all $k$-jets of holomorphic maps 
$X\to Y$. We have $J^0(X,Y)=X\times Y$ and
$$        
	J^1(X,Y)=\{(x,y,\l)\colon x\in X,\ y\in Y,\ 
		\lambda \in {\rm Hom}_{\C} (T_x X,T_y Y)\}.
$$
For a holomorphic map $f\colon X\to Y$ and $k\in \N$ 
we denote by $j^k f\colon X\to J^k(X,Y)$ the $k$-jet extension
of $f$; in particular, $j^0_x f=(x,f(x))$ and 
$j^1_x f = (x,f(x), df_x)$. We shall denote by 
$j^k f|_A$ the restriction of $j^k f$ to points in a 
subset $A\subset X$.

A {\it stratification} of a complex analytic subvariety $A$ 
(in a complex manifold $X$) is a decomposition of $A$ into a 
locally finite disjoint union of open connected complex 
manifolds $A_\a$, called {\it strata} of $A$, 
such that the boundary of each stratum is a union of 
lower dimensional strata (\cite{Wh}, p.\ 227).
Whitney proved  that every complex analyic subvariety
in a complex manifold admits a stratification 
which is {\it regular with respect to tangent planes}  
(\cite{Wh}, theorem 8.5). 
We recall how this {\em Condition (a)} of Whitney
is used in transversality arguments.
Given stratified subvarieties $A\subset X$, $B\subset Y$,
let $\mathcal{NT}_{A,B}  \subset J^1(X,Y)$ 
consist of all $(x,y,\l) \in J^1(X,Y)$ 
such that, if $x$ belongs to a stratum $A_\a$ of $A$
and $y$ belongs to a stratum $B_\b$ of $B$ then 
$$
	\l(T_x A_\a) + T_y B_\b\ne T_y Y.
$$ 
(If $x\not\in A$ or $y\not\in B$ then 
$(x,y,\lambda) \not\in \mathcal{NT}_{A,B}$.)
{\it If the stratifications of $A$ and $B$ satisfy 
Whitney's Condition (a) then $\mathcal{NT}_{A,B}$ is closed
in $J^1(X,Y)$} (\cite{Tr}; \cite{GM}, p.\ 38).
The set $\mathcal{NT}_{A,B}$ is also closed in 
$J^1(X,Y)$ if $B$ is a closed smooth (not necessarily complex) 
submanifold of $Y$; in such case we consider $B$ itself as 
the only stratum.

Let $A\subset X$ and $B\subset Y$ be stratified complex subvarieties. 
Given $f\in \cH(X,Y)$, we say that 
{\it $f|_A$ is transverse to $B$} if the range 
of $j^1 f \colon X\to J^1(X,Y)$ does not intersect 
the set $\mathcal{NT}_{A,B}$. Equivalently, 
$$ 
	(x\in A_\a, \ f(x)\in B_\beta) \,\Longrightarrow\,
	 df_x(T_x A_\a)+ T_{f(x)} (B_\b)=T_{f(x)} Y.
$$

We shall base our discussion on the following condition, 
similar to the one introduced by Gromov (\cite{Gr1}, p.\ 71),  
who indicated its application to transversality theorems
(\cite{Gr1}, p.\ 73, (C')). 

%
%
%
%
\begin{definition}
\label{Ell1} 
Let $X$ and $Y$ be holomorphic (resp.\ algebraic) manifolds.
Holomorphic (resp.\ algebraic) maps $X\to Y$ satisfy 
{\em Condition ${\rm Ell}_1$} if for every holomorphic 
(resp.\ algebraic) map $f\colon X\to Y$ there is a 
holomorphic (resp.\ algebraic) map $F\colon X\times \C^N\to Y$ 
for some $N\ge \dim Y$ such that $F(\cdotp,0)=f$ and 
$F(x,\cdotp)\colon \C^N\to Y$ has rank equal to 
$\dim Y$ at $0\in\C^N$ for every $x\in X$.
\end{definition}

For validity of Condition ${\rm Ell}_1$ see 
proposition \ref{p4.6} below, and also \cite{Gr1}, p.\ 72.

In the sequel, a {\em Whitney stratification} of a complex
analytic subvariety will mean a stratification satisfying 
Whitney's Condition (a).

\begin{theorem} 
\label{t4.2} 
Let $X$ and $Y$ be complex manifolds, with $X$ Stein.
If holomorphic maps $X\to Y$ satisfy Condition ${\rm Ell}_1$ then
for every pair of closed, Whitney stratified complex analytic
subvarieties $A\subset X$, $B\subset J^k(X,Y)$ $(k= 0,1,\ldots)$
the set 
$$
	\{f\in \cH (X,Y)\colon j^k f|_{A} {\rm\ is\ transverse\ to\ } B\}
$$ 
is residual in $\cH(X,Y)$. The same holds if $B$ is a 
smooth closed submanifold of $J^k(X,Y)$. For $k=0$ the 
conclusion holds even if $X$ is not Stein.  
\end{theorem}

In the algebraic category we have the analogous result but
only on compact sets in the source manifold:

\begin{theorem}
\label{t4.3}
Let $X$ and $Y$ be algebraic manifolds, with $X$ affine algebraic.
If algebraic maps $X\to Y$ satisfy Condition ${\rm Ell}_1$ then 
for every compact set $K\subset X$ and Whitney stratified 
complex subvarieties $A\subset X$, 
$B\subset J^k(X,Y)$ the set 
$$
	\{f\in \cO(X,Y)\colon 
        j^k f|_A {\rm \ is\ transverse\ to\ } B {\rm\ on\ } A\cap K \}	
$$
is open and dense in $\cO(X,Y)$. The same holds if $B$ is a 
smooth closed submanifold of  $J^k(X,Y)$. For $k=0$ the 
conclusion holds without assuming that $X$ be affine. 
\end{theorem}

Theorems \ref{t1.3} and \ref{t1.4} in \S 1 follow immediately by 
combining theorems \ref{t4.2}, \ref{t4.3} and proposition \ref{p4.6} below.

\begin{proof} 
The proofs of theorems \ref{t4.2} and \ref{t4.3} 
are parallel up to the point where the 
Baire property of the space of $\cH(X,Y)$ is invoked; 
in the algebraic case this leaves us with the weaker statement. 
We shall follow Abraham's reduction \cite{Ab} to Sard's theorem 
\cite{Sa}; although this is well known (see \S 1.3.7.\ in \cite{GM} 
as well as \cite{Fo}, \cite{KZ}), we include a sketch of proof 
since we shall need small modifications of the standard arguments 
in some of the proofs in this section and in \S 5.

Whitney's Condition (a) implies the following. 

%
%
\begin{lemma} 
\label{l4.4} 
Let $X$ and $Y$ be complex manifolds and let $A\subset X$,
$B\subset J^k(X,Y)$ be closed, Whitney stratified complex 
analytic subvarieties. For every compact subset  $K$ of $X$ the set 
$$
   \cT_{A,B,K}= \{f\in \cH (X,Y) \colon j^k f|_{A}\ 
        {\rm\ is\ transverse\ to\ } B {\rm\ on\ } A\cap K \}
$$
is open in $\cH(X,Y)$.  The analogous result holds if $B$ 
is a smooth closed submanifold.  
\end{lemma}

\begin{proof}
Consider the basic case with $B\subset Y$. 
Given a map $f\colon X\to Y$ and a compact set $K\subset X$, 
the restriction $f|_A\colon A\to B$ is transverse to $B$ 
at each point of $A\cap K$ if and only if 
$(j^1 f)(K) \cap \mathcal{NT}_{A,B} =\emptyset$.
Assuming this to be the case, and taking into account that 
$\mathcal{NT}_{A,B}$ is closed in $J^k (X,Y)$ by Whitney's condition, 
there is a compact set $L\subset X$, with $K\subset {\rm Int}\, L$,
such that $(j^1 f)(L)\cap \mathcal{NT}_{A,B} =\emptyset$. 
If $g \in \cH(X,Y)$ is sufficiently uniformly close to $f$ on $L$ 
then $j^1 g$ is close to $j^1 f$ on $K$
and hence $(j^1 g)(K)\cap \mathcal{NT}_{A,B}=\emptyset$. 
In the general case one applies the same 
argument with $f$ replaced by the map 
$j^k f\colon X\to \wt Y=J^k(X,Y)$.
\end{proof}

To prove theorem \ref{t4.2} it suffices to show that for 
every compact $K$ in $X$ the set $\cT_{A,B,K} \subset \cH (X,Y)$ 
(which is open in $\cH (X,Y)$ by lemma \ref{l4.4}) is 
everywhere dense in $\cH(X,Y)$. Since $\cH(X,Y)$ is a Baire space, 
the conclusion of theorem \ref{t4.2} then follows by taking 
the intersection of such sets over a countable family 
of compacts exhausting $X$. In the algebraic case 
we omit the last step. 

Consider first the basic case with $A=X$ and $B\subset Y$.
Let $f\colon X\to Y$ be a holomorphic (resp.\ algebraic) map.
Choose a map $F\colon X\times \C^N\to Y$ as in definition 
\ref{Ell1} of Condition ${\rm Ell}_1$.
Let $\pi\colon X\times \C^N\to \C^N$ denote the projection 
$\pi(x,t)=t$.  Fix a compact set $K$ in $X$. Since 
$\partial_t F(x,0) \colon T_0\C^N \to T_{f(x)} Y$ 
is surjective for every $x\in X$, there are a small
ball $D\subset \C^N$ around the origin and an open 
set $U\subset X$ containing $K$ such that $F$ is a submersion 
of $V= U\times D$ to $Y$. Hence $B'= F^{-1}(B)\cap V$ 
is a closed, Whitney stratified, complex analytic subvariety 
of $V$ (the strata $B_\b$ of $B$ pull back by $F|_V$ to the 
strata $B'_\b$ of $B'$). 
Set $f_t= F(\cdotp,t)\colon X\to Y$ for $t\in \C^N$. 
If $(x,t) \in B'_\b$ then $y=f_t(x) \in B_\b$, 
and by inspecting the definitions we see that the 
following are equivalent (compare \cite{GM}, p.\ 40): 

\begin{itemize}
\item[(a)] $(df_t)_x (T_x X) + T_{y} B_\b = T_{y} Y$;
\item[(b)] $(x,t)$ is a regular point of the
restricted projection $\pi \colon B'_\b \to D$. 
\end{itemize}

By Sard's theorem \cite{Sa}, applied inductively to the components
of a projection $\pi$ in (b), we see that the set of regular values of all
projections in (b) is residual in $D$.  
Choosing $t$ in this set and close to $0$ 
we get a map $f_t\colon X\to Y $ which is transverse to $B$ on $U$
and which approximates $f=f_0$ uniformly on $K$. The same argument
applies if $B$ is a smooth closed submanifold of $Y$.

If $A$ is a Whitney stratified complex subvariety of $X$, 
one applies the above argument with $U$ replaced 
by $U\cap A_\a$ for a fixed stratum $A_\a$ of $A$
($f$ and $F$ are still defined globally on $X$). 
This gives a residual set of $t$'s in $D\subset \C^N$ 
for which $f_t|_{A_\a \cap U}$ is transverse to a stratum $B_\b$. 
Since $A$ and $B$ have at most countably many strata
and $\C^N$ is a Baire space, we find $t\in \C^N$ arbitrarily 
close to $0$ such that $f_t|_{A\cap U}$ is transverse to $B$. 

This proves the basic transversality theorem for holomorphic maps 
$X\to Y$.  All steps hold for algebraic maps as well, even if 
the subvarieties $A$ and $B$ are non-algebraic, and we did not
need any special properties of $X$ and $Y$ other than Condition 
${\rm Ell}_1$ for algebraic maps $X\to Y$.   

To prove theorems \ref{t4.2} (resp.\ \ref{t4.3}) for $k=0$
we must consider transversality of holomorphic 
(resp.\ algebraic) sections $x\to j^1_x f = (x,f(x))$ 
$(x\in X)$ of the trivial fibration $X\times Y\to X$ 
with respect to complex (resp.\ algebraic) subvarieties 
of $X\times Y$. This is done by obvious modifications 
of the above arguments, using the fact that the map 
$(x,t)\to (x,F(x,t)) \in X\times Y$, constructed above, 
is a submersion on the subset $U\times D\subset X\times \C^N$.

Consider now the case $k>0$.
Fix a map $f\colon X\to Y$ and a compact set $K\subset X$.
The goal is to prove that $f$ can be approximated uniformly
on $K$ by holomorphic (resp.\ algebraic) maps $X\to Y$ 
whose $k$-jet extension $j^k f$ is such that 
$j^k f|_A \colon A\to J^k(X,Y)$ is transverse to 
the subvariety $B\subset J^k(X,Y)$ at each point of $A\cap K$. 

Let $F\colon X\times\C^N\to Y$ be as in definition \ref{Ell1},
with $F(\cdotp,0)=f$. Recall that $X$ is assumed to be Stein
(in theorem \ref{t4.2}) resp.\ affine algebraic
(in theorem \ref{t4.3}). Thus we may assume that $X$ is a closed
holomorphic (resp.\ algebraic) submanifold of a Euclidean space $\C^n$. 
Let $\cW$ denote the complex vector space of all polynomial maps 
$P\colon \C^n\to\C^N$ of degree $\le k$. Consider the holomorphic 
(resp.\ algebraic) map $G \colon X\times \cW\to Y$ defined by
$$
	G(x,P)=F(x,P(x)), \qquad x\in X,\ P\in \cW.
$$
For each $P\in \cW $ set $G_P=G(\cdotp, P)\colon X\to Y$;
then $G_0(x) = F(x,0)=f(x)$.

\begin{lemma} 
\label{l4.5} 
The map $\Phi \colon  X\times \cW  \to J^k(X,Y)$, defined by 
$\Phi(x,P)= j^k_x (G_P)$, is a submersion in an open 
neighborhood of $X\times \{0\}$ in $X\times \cW$. 
\end{lemma}

\begin{proof} 
The argument is local and hence we may assume that $X=\C^n$. 
Write $P=(P_1,\ldots,P_N)\in \cW$, and let $t=(t_1,\ldots,t_N)$ 
be coordinates on $\C^N$. For every multiindex 
$I=(i_1,\ldots,i_n)$ we have 
$$ 
	\partial^I_x (G_P) = \sum_{j=1}^N 
	\frac{\partial}{\partial t_j} F(x,P(x))\, \partial^I_x P_j(x) + H_I(x),
$$
where $H_I(x)$ contains only terms $\partial^J_x P$, with $|J|<|I|$,
multiplied by various partial derivatives of $F$.
Hence the $k$-jet map $j^k_x (G_P)$ is lower triangular
with respect to the components of $j^k_x P$, and the diagonal 
terms are nondegenerate at $P=0$ (since $G_0(x)=F(x,0)$ and 
$\partial_t F(x,0)$ is nondenegenerate), thus proving the lemma. 
\end{proof}

Sard's theorem, applied to the map $\Phi$ in lemma \ref{l4.5},
shows that for most $P\in \cW$ the map $j^k (G_P)|_A$ is transverse to 
the subvariety $B\subset J^k(X,Y)$ at every point of $A\cap K$. 
This proves theorems \ref{t4.2} and \ref{t4.3}.
\end{proof}

Combining theorems \ref{t4.2} and \ref{t4.3} with the following
proposition gives several  transversality theorems, 
including those announced in \S 1. Compare with the
examples in \cite{Gr1}, p.\ 72.

\begin{proposition} 
\label{p4.6} 
Let $X$ and $Y$ be complex manifolds.
\begin{enumerate}
\item [(a)]
If $Y$ admits a dominating spray $s\colon Y\times \C^N\to Y$ 
defined on a trivial bundle then holomorphic maps 
$X\to Y$ satisfy Condition ${\rm Ell}_1$. This holds
if $Y$ is a complex homogeneous space. 
\item [(b)]
Let $X$ be Stein. If $Y$ is subelliptic (definition \ref{def1}) or, 
more generally, if $Y$ enjoys {\rm CAP} then holomorphic maps 
$X\to Y$ satisfy Condition ${\rm Ell}_1$. 
\item [(c)]
If $X$ and $Y$ are algebraic and if $Y$ admits a dominating 
algebraic spray $s\colon Y\times \C^N\to Y$ then algebraic maps 
$X\to Y$ satisfy Condition ${\rm Ell}_1$.
\item [(d)]
If $X$ is affine algebraic and $Y$ is algebraically subelliptic 
then algebraic maps $X\to Y$ satisfy Condition ${\rm Ell}_1$.
\end{enumerate}
\end{proposition}

\begin{proof} 
Fix a holomorphic map $f\colon X\to Y$. If $Y$ admits a dominating 
spray $(E,p,s)$ then $f^* E\to X$ is a holomorphic vector bundle, 
and there is a fiberwise bijective holomorphic map 
$\iota\colon f^* E\to E$ covering $f$.  
The map  $F=s\circ\iota \colon f^* E\to Y$ satisfies 
Condition ${\rm Ell}_1$ for $f$, except that $f^* E$ need 
not be a trivial bundle over $X$. In case (a) the bundle $E\to Y$
is assumed to be trivial, hence $f^*E$ is also trivial and 
(a) follows. If $Y$ is a homogeneous space of a complex Lie group 
$G$, with Lie algebra ${\mathbf g} =T_e G$, the map
$s\colon Y\times {\bf g}\to Y$, $s(y,v)=\exp(v)\cdotp y$,
is a dominating spray defined on a trivial bundle.
The analogous argument proves (c). 

If $X$ is Stein then by Cartan's Theorem A 
the holomorphic vector bundle $f^*E\to X$ 
is generated  by finitely many (say $N$) holomorphic sections, 
and hence there is a surjective complex vector bundle map 
$\tau \colon X\times \C^N\to f^*E$. 
The map $F= s\circ\iota \circ \tau \colon X\times \C^N\to Y$ 
satisfies Condition ${\rm Ell}_1$ with respect to 
$f=F(\cdotp,0)$. The analogous argument holds in the algebraic case
by appealing to Serre's Theorem A (\cite{Se1}, p.\ 237, Th\'eor\`eme 2),
thus proving (b) (resp.\ (d)) for (algebraically) 
elliptic target manifolds.

Assume now that $Y$ is (algebraically) subelliptic 
and let $X$ be an affine manifold (Stein resp.\ affine algebraic). 
Let $(E_j,p_j,s_j)$ for $j=1,\ldots, k$ be a finite dominating 
family of holomorphic (resp.\ algebraic) sprays on $Y$
(def.\ \ref{def1}). We shall perform the above procedure several times. 
In essence we use dominating composed sprays as in the 
proof of theorem \ref{t3.3} above; see also \cite{Gr2}, \S 1.3.A', and 
\cite{FP1}, corollary 3.5. 

Let $E'_1 = f^* E_1 \to X$ be the pull-back of 
$\pi_1\colon E_1\to Y$ by $f\colon X\to Y$, and define 
$\sigma_1\colon E'_1\to Y$ by $\sigma_1(x,e)=s_1(f(x),e)$.  
As before, there is a surjective complex vector bundle map 
$X\times \C^{n_1}\to E'_1$ for some $n_1\in \N$.
By composing it with $\sigma_1$ we obtain a map
$f_1\colon X_1 = X\times \C^{n_1}\to Y$
satisfying $f_1(x,0)=f(x)=y\in Y$ and 
$$
	\partial_t f_1(x,t)|_{t=0} (T_0 \C^{n_1}) = (ds_1)_y (E_{1,y}) 
	\subset T_y Y.
$$
Repeating the construction with $f_1\colon X_1 \to Y$ 
and the second spray $s_2\colon E_2\to Y$ 
we find an integer $n_2\in \N$ and a holomorphic map 
$f_2\colon X_2 = X_1\times \C^{n_2}= X\times \C^{n_1}\times\C^{n_2}\to Y$
satisfying $f_2(x,t,0)=f_1(x,t)$ (hence $f_2(x,0,0)=f(x)=y$) and 
$$
	\partial_u f_2(x,0,u)|_{u=0} (T_0 \C^{n_2}) = 
	(ds_2)_y (E_{2,y}) \subset T_y Y.
$$
After $k$ steps we obtain a map $F\colon X\times \C^N\to Y$
$(N=n_1+\cdots + n_k)$ satisfying the following for every 
$x\in X$ and $y=f(x)\in Y$:
$$
	F(x,0)=f(x), \quad     
	\partial_t F(x,t)|_{t=0} (T_0\C^N) = 
	\sum_{j=1}^k (ds_j)_{y}(E_{j,y}) = T_y Y.
$$
The last equality is the domination property (2.1).
This completes the proof of (b) for a subelliptic $Y$;
the same proof applies in the algebraic case (d)
by appealing to the Theorem A of Serre \cite{Se1} when 
passing at each step to a trivial bundle.

It remains to prove (b) when $X$ is Stein and $Y$ 
enjoys CAP (which is equivalent to the Oka property
by corollary \ref{hierarchy}). Let $f\colon X\to Y$ be a holomorphic
map. Consider the associated embedding $x\in X \to (x,f(x))\in X\times Y$
with normal bundle $E=f^* TY\to X$. 
By the Docquier-Grauert theorem \cite{DG} there are 
an open neighborhood $V\subset E$ of the zero section $X\subset E$
and a biholomorphic map $G \colon V\to G(V)\subset X\times Y$
of the form $G(x,\xi)=(x,g(x,\xi))$ $(x\in X,\ \xi \in E_x)$ 
satisfying $g(x,0_x)=f(x)$. We can extend $g$ to a continuous
map $E\to Y$ without changing its values on a smaller 
neighborhood of the zero section $X\subset E$. 

Since $Y$ enjoys CAP, it also enjoys the Oka property 
with jet interpolation (corollary 1.4 in \cite{F7};
see also remark \ref{generalOka} above).
As $E$ is a Stein manifold, this gives a holomorphic map 
$\wt g \colon E\to Y$ which agrees with $g$ to the second
order along the zero section $X\subset E$. 
Let $\iota\colon X\times\C^N\to E$ be a surjective
holomorphic vector bundle map (which exists for large $N$
by Cartan's Theorem A). The composition 
$F=\wt g\circ\iota \colon X\times \C^N\to Y$ then
satisfies definition \ref{Ell1} with respect to the map $f$.
\end{proof}

\begin{corollary}
\label{c4.7}
Holomorphic maps from any complex manifold $X$ to a complex 
homogeneous manifold $Y$ satisfy the basic transversality theorem 
(for zero-jets). The same holds if $Y=\C^n\bs A$ where $A$ is a thin
algebraic subvariety, and in this case the basic transversality theorem 
also holds for algebraic maps from any algebraic manifold to $\C^n\bs A$.
\end{corollary}

In the case $Y=\P_n$ we recover Bertini's theorem (\cite{GM}, p.\ 150).

%
%
%
%
In \cite{KZ} Kaliman and Zaidenberg proved the following 
theorem in which there is no restriction on the target manifold, 
but the domain of the map shrinks.

\begin{theorem}
\label{t4.8} {\rm (\cite{KZ})} 
Assume that $X$ is a Stein manifold, $Y$ is a complex manifold 
and $A\subset X$, $B\subset J^k(X,Y)$ are 
closed, Whitney stratified complex subvarieties. 
For any $f\in \mathcal{H}(X,Y)$ and any compact set $K\subset X$ 
there is a holomorphic map $g\colon U\to Y$ in an open 
neighborhood of $K$ such that $j^k g|_{A\cap U}$ 
is transverse to $B$, and $g$ approximates $f$ as close 
as desired uniformly on $K$. 
\end{theorem}

Theorem \ref{t4.8} also follows from our proof of theorem \ref{t4.2}
by using a holomorphic map $F\colon U\times D\to Y$ satisfying 
Condition ${\rm Ell}_1$ along $U\times \{0\}$, where 
$U\subset X$ is an open neighborhood of $K$ 
and $D\subset \C^N$ is a small ball around $0\in\C^N$. 
Such $F$ exists provided that $K$ admits a basis of
open Stein neighborhoods in $X$ (which is  the case
if $K$ is $\cH(X)$-convex). 
Indeed, the set $\{(x,f(x)) \colon x\in K\} \subset X\times Y$
has an open Stein neighborhood $\Omega\subset X\times Y$ \cite{Si1}, \cite{D1},
and hence there exist holomorphic vector fields $V^1,\ldots,V^N$ 
in $\Omega$ tangent to the fibers of the projection 
$X\times Y \to X$ and generating the tangent 
space of the fiber at each point. Let $\theta^j_t$ denote 
the flow of $V^j$. For $x\in X$ near $K$ 
and small $t_1,\ldots,t_N\in \C$ the map 
$ 
     F(x,t_1,\ldots,t_N)= 
     \pi_Y \circ \theta^1_{t_1}\circ\cdots\circ\theta^N_{t_N}(x,f(x))
$
satisfies the required property.

\medskip
\noindent\textit{Alternative proof of theorem \ref{t1.4}}.
Assume that $X$ is Stein and $Y$ enjoys the Oka property. 
Let $f\colon X\to Y$ be a holomorphic map. 
Choose compact $\cH (X)$-convex subsets $K,L \subset X$
with $K\subset {\rm Int}\, L$. By theorem \ref{t4.8} we 
can approximate $f$ uniformly on $L$ by a holomorphic map
$g\colon U\to Y$ on an open set $U\supset L$
such that $j^k g|_{U\cap A}$ is transverse to $B$. 
If the approximation is sufficiently close, there is a smooth  
map $\widetilde g\colon X\to Y$ which agrees with $g$ in a 
neighborhood of $L$, and agrees with $f$ on $X\backslash U$. 
By the Oka property of $Y$ the map 
$\widetilde g$ can be approximated uniformly on $L$ 
by holomorphic maps $\widetilde f\colon X\to Y$. If the approximation 
is sufficiently close then $\widetilde f$ still satisfies the 
desired transversality condition on $K$ by lemma \ref{l4.4}. 
This shows the density of transverse maps on compacts in $X$, 
thus completing the proof.
\smallskip

We now give an interpolation version of theorem \ref{t1.4}. 
Given a closed complex subvariety $X_0$ of $X$, 
$f_0\in \cH (X,Y)$ and $r\in \{0,1,\ldots\}$, the set
$$
	\cH (X,Y;X_0,f_0,r) = 
	\{f\in \cH (X,Y)\colon j^r f|_{X_0}= j^r f_0|_{X_0}\}
$$
is a closed metric subspace of $\cH (X,Y)$, hence a Baire space.

\begin{theorem}
\label{t4.9} 
Let $X$ be a Stein manifold and $Y$ a complex manifold 
enjoying the Oka property. 
Let $A\subset X$ and $B\subset J^k(X,Y)$ be closed, Whitney 
stratified complex subvarieties. If $f_0\in \cH (X,Y)$ 
is such that $j^k f_0|_{A}$ is transverse to $B$ at all points of 
$A\cap X_0$ then for every integer $r\ge k$ there is a residual set 
of $f\in \cH (X,Y;X_0,f_0,r)$ for which $j^k f|_{A}$ is 
transverse to $B$.  
\end{theorem}

\begin{proof} 
Since $r\ge k$, the set of all $f\in \cH (X,Y;X_0,f_0,r)$
for which $j^k f|_{A}$ is transverse to $B$ at all points of 
$A\cap X_0$ is open in $ \cH (X,Y;X_0,f_0,r)$.
To prove theorem \ref{t4.9} it thus suffices to show 
that we can approximate the initial map $f_0$ uniformly on 
any compact $\cH(X)$-convex subset $K\subset X$  
by $f \in \cH (X,Y;X_0,f_0,r)$ such that $j^k f|_{A}$ 
is transverse to $B$ at every point of $A\cap K$. 

Assume first that holomorphic maps $X\to Y$ satisfy 
Condition ${\rm Ell}_1$, and let 
$F\colon X\times \C^N\to Y$ be as in definition
\ref{Ell1}, with $F(\cdotp,0)=f_0$. 
Consider the basic case $k=0$, $B\subset Y$. 
There exist functions $g_1,\ldots,g_l\in \cH (X)$ 
which vanish to order $r+1$ on the subvariety 
$X_0=\{x\in X\colon g_j(x)=0,\ j=1,\ldots,l\}$. 
For every $x\in X$ let $\sigma_x \colon (\C^N)^l \to \C^N$ be defined by 
$$ 
	\sigma_x(t_1,\ldots,t_l) = \sum_{j=1}^l t_j g_j(x), \qquad
	t_j\in \C^N,\ j=1,2,\ldots, l.
$$
Clearly $\sigma_x$ is surjective if $x\in X\backslash X_0$ and 
is the zero map if $x\in X_0$. 
The map $\widetilde F\colon X\times \C^{Nl} \to Y$, defined by
$\widetilde F(x,t)= \widetilde F(x,t_1,\ldots,t_l)=F(x,\sigma_x(t_1,\ldots,t_l))$,
is a submersion with respect to $t$ (at $t=0$) if $x\in X\backslash X_0$, and is 
degenerate (constant) if $x\in X_0$. 
Hence the proof of theorem \ref{t4.2} applies 
over $X\backslash X_0$. 

Let $f_t =\widetilde F(\cdotp, t) \colon X\to Y$ for $t\in \C^{Nl}$.
By construction $j^r f_t|_{X_0}= j^r f_0|_{X_0}$ for every $t$.
Choose a compact set $K\subset X$. By the assumption 
$f_0|_{A}$ is transverse to $B$ on $A\cap X_0$.
Hence there is an open neighborhood $U\subset X$ of $A\cap X_0\cap K$
such that $f_t|_{A\cap U}$ is transverse to $B$ for every $t$ 
sufficiently close to $0$ (lemma \ref{l4.4}).
The set $K'=K\backslash U \subset X\backslash X_0$ 
is compact and hence for most values of $t$ the map $f_t|_A$ 
is transverse to $B$ on $A\cap K'$. Thus $f_t|_A$ is transverse 
to $B$ on $K\cap A$ for most $t$ close to $0$ which concludes
the proof for $k=0$. Similarly one obtains the proof for $k>0$
by following the arguments in the proof of  theorem \ref{t4.2}.

The same proof gives a semiglobal version of theorem \ref{t4.9},
analogous to theorem \ref{t4.8}, without any restriction on 
the manifold $Y$. 

The proof of the general case is completed exactly 
as the alternative proof of theorem \ref{t1.4} given above,
using corollary 1.4 in \cite{F7} to the effect that the 
(basic) Oka property of $Y$ implies the Oka property with 
jet interpolation on a closed complex subvariety $X_0$
of a Stein manifold $X$.
\end{proof}

In the algebraic category the global transversality theorem
holds under the following stronger assumption on $Y$.

\begin{proposition}
\label{p4.10} 
If $Y$ is an algebraic manifold with a submersive algebraic spray
$s\colon E\to Y$ (i.e., such that $s\colon E_y\to Y$ is a submersion 
for every $y\in Y$) then algebraic maps $X\to Y$ from any 
affine algebraic manifold $X$  satisfy the jet transversality 
theorem with respect to closed complex analytic subvarieties.
\end{proposition}

\begin{proof} 
Let $f_0\colon X\to Y$ be an algebraic map.
Pulling back the submersive algebraic spray $s\colon E\to Y$ by $f_0$ 
we obtain an {\it algebraic submersion} $F \colon X\times \C^N\to Y$
satisfying $f_0=F(\cdotp,0)$ (compare with the proof of proposition 
\ref{p4.6}). Given closed complex subvarieties 
$A\subset X$ and $B\subset Y$, Sard's theorem
shows that for a generic choice of $t\in\C^N$
the algebraic map $f_t|_A = F(\cdotp, t)|_A$ is 
transverse to $B$ (see the proof of theorem \ref{t4.2}). 
Similarly one obtains the jet transversality theorem by considering 
maps $x\to F(x,P(x))$ for polynomials $P \colon \C^n\to \C^N$,
where $X$ is embedded in $\C^n$.  
\end{proof}

%
%
%
%
%
\section{The homotopy principle for holomorphic submersions} 
Given complex manifolds $X$ and $Y$, we denote by $\cS (X,Y)$ 
the set of all pairs $(f,\iota)$ where $f\colon X\to Y$ is a continuous 
map and $\iota\colon TX\to TY$ is a fiberwise surjective complex vector 
bundle map such that the diagram in Fig.\ 1 commutes.  

\begin{figure}[ht]
$$ 
  \begin{array}[]{*3c}  
            TX         & \stackrel{\iota}{\longrightarrow} & TY        \\
            \downarrow{} &                  & \downarrow{} \\
             X         & \stackrel{f}{\longrightarrow}   & Y    
  \end{array}
$$
\centerline{Figure 1: The space $\cS (X,Y)$} 
\end{figure}

The existence of $\iota$ in Fig.\ 1 (for a given $f$) 
is invariant under homotopies of the base map, and 
$Tf = (f,df) \in \cS (X,Y)$ precisely when $f\colon X\to Y$ is a 
holomorphic submersion. (Here $Tf$ denotes the tangent map of $f$.)
Hence a necessary condition for a continuous map $f\colon X\to Y$ 
to be homotopic to a holomorphic submersion $X\to Y$ is 
that $f$ be covered by a map $\iota\colon TX \to TY$ 
such that $(f,\iota)\in \cS (X,Y)$. 
By theorem II in \cite{F3} this condition is also sufficient
if $X$ is Stein and $Y=\C^p$, $p<\dim X$. Here we prove
the same result when the target manifold is of Class $\cA$,
or a holomorphic quotient of such manifold. 

\begin{theorem} 
\label{t5.1}  
Assume that $Y$ is a complex manifold which admits an
unramified holomorphic covering $\widetilde Y\to Y$ 
by a quasi-projective algebraic manifold $\wt Y$ 
of Class $\cA$ (def.\ \ref{def2}). If $X$ is a 
Stein manifold with $\dim X > \dim Y$ and $K\subset X$ 
is a compact $\cH (X)$-convex subset then the following
hold.
\begin{itemize}
\item[(a)] 
For any $(f,\iota)\in \cS (X,Y)$ the map $f$ is 
homotopic to a holomorphic submersion $f_1\colon X\to Y$.
If in addition $f|_K\colon K\to Y$ is a holomorphic submersion
and $\iota|_K=df|_K$ then $f_1$ can be chosen to
approximate $f$ uniformly on $K$. 
\item[(b)] 
Holomorphic submersions $f_0,\, f_1\colon X\to Y$
are regularly homotopic through a family of holomorphic submersions 
$X\to Y$ if and only if their tangent maps $Tf_0$ and $Tf_1$ 
belong to the same path connected component of $\cS(X,Y)$.
\item[(c)] If $\dim X \ge 2\, \dim Y-1$ then every continuous map 
$X\to Y$ is homotopic to a holomorphic submersion;  
if $\dim X\ge 2\, \dim Y$ then any two holomorphic 
submersions $X\to Y$ are regularly homotopic.
\end{itemize}
\end{theorem}

Theorem \ref{t5.1} implies results on the existence of nonsingular
holomorphic foliations of the source manifold $X$; 
see \cite{F3} for the case $Y=\C^p$.

\begin{example} 
We list some examples of (holomorphic quotients of) 
quasi-pro\-jec\-tive algebraic manifolds of Class $\cA$
to which theorem \ref{t5.1} applies.
\begin{enumerate}
\item 
$Y=\widehat Y\bs A$ where $\widehat Y$ is an affine space,
a projective space or a Grassmanian and $A$ is a 
thin algebraic subvariety of $\widehat Y$.
\item
$Y=\C^p/\Gamma$ where $\Gamma$ is a lattice in $\C^p$. 
This class includes all complex tori.
\item
{\em Hopf manifolds} are quotients of 
$\C^p_* \stackrel{\rm def}{=} \C^p\backslash\{0\}$ $(p\ge 2)$ by an 
infinite cyclic group, or a finite extension of such 
group (\cite{BH}, p.\ 225).
\item
Let $\pi\colon W\to Y$ be a holomorphic fiber bundle 
whose base $Y$ is a quotient of a Class $\cA _0$ manifold, 
the fiber $\pi^{-1}(y)$ is $\C^m$ respectively $\P_m$, 
and the structure group is $GL_m(\C)$ respectively $PGL_m(\C)$. 
It is easily seen that $W$ is then a quotient of a Class $\cA_0$ manifold. 
\end{enumerate}
\end{example}

Part (c) of theorem \ref{t5.1} follows from (a) and (b) by 
topological reasons (see corollary 2.3 in \cite{F4}). 
We shall reduce parts (a) and (b) to theorem 2.1 in \cite{F4}. 
To this end we must recall from \cite{F4} a certain holomorphic
flexibility property, called Property ${\mathrm S}_{\mathrm n}$,
which is the localization to Euclidean spaces of the 
homotopy principle for holomorphic submersions from $n$-dimensional 
Stein manifolds to $Y$. 

Let $z=(z_1,\ldots,z_n)$, $z_j=x_j+ iy_j$, 
denote the coordinates on $\C^n$. Set 
$$
	P =\{z\in \C^n\colon |x_j| \le 1,\ |y_j|\le 1,\ j=1,\ldots,n\}.  \eqno(5.1)
$$
A compact convex subset $K\subset \C^n$ is {\em special} if
$$ 
	K=\{z\in P \colon y_n \le h(z_1,\ldots,z_{n-1},x_n)\}   \eqno(5.2)
$$
where $h$ is a smooth concave function with  values in $(-1,+1)$.

%
%
\begin{definition}
\label{def3}
Let $Y$ be a complex manifold and $d$ a distance function 
on $Y$ induced by a Riemannian metric. 
\begin{itemize}
\item[(a)]
$Y$ satisfies {\em Property ${\mathrm S}_{\mathrm n}$} 
if for every holomorphic submersion 
$f\colon K\to Y$ on a special compact convex set $K$ (5.2) and for every 
$\e>0$ there is a  holomorphic submersion $\widetilde f\colon P\to Y$
satisfying $\sup_{x\in K} d(f(x),\widetilde f(x)) <\e$.
\item[(b)] 
$Y$ satisfies {\em Property $\mathrm{HS}_{\mathrm n}$} if for every homotopy of 
holomorphic submersions $f_t\colon K\to Y$ $(t\in [0,1])$ 
such that $f_0$ and $f_1$ extend to holomorphic submersions 
$P\to Y$ there is for every $\e>0$ a homotopy of holomorphic 
submersions $\widetilde f_t\colon P \to Y$ $(t\in [0,1])$ satisfying
$\widetilde f_0=f_0$, $\widetilde f_1=f_1$, and  
$\sup\{ d(f_t(x),\widetilde f_t(x))\colon x\in K,\ t\in [0,1]\} <\e$.
\end{itemize}
\end{definition}

According to theorem 2.1 in \cite{F4}, the conclusion of (a) 
(respectively of (b)) in theorem \ref{t5.1} holds provided 
that $Y$ satisfies Property ${\mathrm S}_{\mathrm n}$ 
(respectively Property $\mathrm{HS}_{\mathrm n}$) 
with $n=\dim X$. Furthermore, both properties are 
obviously invariant when passing to an unramified holomorphic 
covering or quotient. Hence theorem \ref{t5.1} is implied 
by the following proposition.

\begin{proposition} 
\label{p5.4}
A quasi-projective algebraic manifold $Y$ of Class $\cA$ 
(def.\ \ref{def2}) satisfies {\em Properties} 
${\mathrm S}_{\mathrm n}$ and $\mathrm{HS}_{\mathrm n}$ 
for every $n> \dim Y$.
\end{proposition}

\begin{proof} 
Let $Y$ be a manifold of class $\cA$.
We may (and shall) assume that $Y$ is connected.
We first prove Property ${\mathrm S}_{\mathrm n}$. 
Let $f\colon K\to Y$ be a holomorphic submersion from an open 
neighborhood of a special compact convex $K\subset \C^n$
(5.2). By corollary \ref{c2.4} the manifold 
$Y$ is algebraically subelliptic, and corollary \ref{c3.2}
then shows that $f$ can be approximated uniformly on $K$
by algebraic maps $\C^n\to Y$. 
Thus it suffices to consider the case when $f\colon \C^n\to Y$ is 
an algebraic map which is a submersion at every point of $K$.

We denote by $\Sigma\subset \C^n$ the ramification locus of $f$
(the set of nonsubmersion points). This is an algebraic 
subvariety of $\C^n$ which does not intersect $K$.

\smallskip
{\em Case 1:} $\dim \Sigma \le n-2$.  
Lemma 3.4 in \cite{F3} provides a holomorphic automorphism $\psi$ 
of $\C^n$ which approximates the identity map in a neighborhood of $K$ 
and satisfies $\psi(P) \cap \Sigma =\emptyset$. 
The map $\widetilde f=f \circ \psi \colon P\to Y$
is a holomorphic submersion approximating $f$ on $K$,
thus proving Property ${\mathrm S}_{\mathrm n}$.
For this argument it suffices to assume that $f$ is a 
holomorphic submersion (with values in $Y$) 
defined in the complement of a thin algebraic subvariety 
in $\C^n$.

\smallskip
{\em Case 2:}  \rm $\dim \Sigma=n-1$.
We shall reduce to Case 1 by inductively removing all
$(n-1)$-dimensional irreducible components from 
the ramification locus $\Sigma$, changing the map
at every step. 

Choose an irreducible component $\Sigma' \subset \Sigma$ 
of dimension $n-1$ and a point $z_0\in \Sigma'$ which 
does not belong to any other irreducible component of $\Sigma$. 
By definition \ref{def2} we have $Y=\widehat Y\backslash A$ 
where $\widehat Y$ is a connected manifold of Class $\cA_0$ and $A$ 
is a thin algebraic subvariety of $\widehat Y$.
Let $U\subset \widehat Y$ be a Zariski open set 
isomorphic to $\C^p$ and containing the point $y_0=f(z_0)$. 
Let $s_0\colon U\times \C^p \to U\simeq\C^p$ 
denote the spray $s_0(y,t)=y+t$.  
Choose a closed algebraic subvariety $Y_0\subset \widehat Y$ 
of pure dimension $p-1$ such that 
$\widehat Y= Y_0\cup U$ and $y_0\notin Y_0$.
Let $L=[Y_0]^{-1}$ where $[Y_0] \to \widehat Y$ is the 
holomorphic line bundle defined by the divisor of $Y_0$. 
Let $\tau_p=\widehat Y\times\C^p$. 

By proposition 1.3 in \cite{F2}, p.\ 541
(or \S 3.5.B.\ and \S 3.5.C.\ in \cite{Gr2})
there are an integer $m\in \N$ and an algebraic spray 
$s\colon E = \tau_p\otimes L^{\otimes m} \to \widehat Y$ 
such that $s(y,t)=y$ for all $y\in Y_0$ and $t\in E_y$,
and $s$ equals $s_0$ over the open set 
$\widehat Y\backslash Y_0 \subset U$, using an
identification of $E|_U$ with $\tau_p|_U$. 
(The line bundle $L$ is trivial over $\widehat Y\bs Y_0$.)

By Serre's Theorem A \cite{Se1} the algebraic vector bundle $f^*(E) \to \C^n$ 
is generated by finitely many (say $q$) algebraic sections, 
and hence there is a surjective algebraic vector bundle map 
$\rho\colon \C^n\times \C^q\to f^*E$. Let $\iota\colon f^*E\to E$
be the natural map covering $f$. Set $Z=f^{-1}(Y_0)\subset\C^n$.
The algebraic map 
$F=s\circ\iota\circ\rho \colon \C^n\times \C^q\to \widehat Y$
satisfies the following properties:
\begin{itemize}
\item[(a)] $F(z,0)=f(z)$ $(z\in \C^n)$, 
\item[(b)] $F(z,t)=f(z)$ $(z\in Z\subset \C^n,\ t\in \C^q)$, and
\item[(c)] $F(z,\cdotp)\colon \C^q \to \widehat Y$ is a submersion
for every point $z\in V:=\C^n\backslash Z$.
\end{itemize}

The proof of theorem \ref{t4.2} (see especially lemma \ref{l4.5}
and proposition \ref{p4.10}) 
gives  a polynomial map $P\colon \C^n\to \C^q$
such that the algebraic map $f_1\colon \C^n\to \widehat Y$ 
defined by $f_1(x)=F(x,P(x))$ $(x\in \C^n)$ 
satisfies the following properties:
\begin{itemize}
\item[(i)]   $f_1$ approximates $f$ uniformly on $K$ as close as desired,
\item[(ii)]  $j^1 f|_{Z}= j^1 f_1|_{Z}$,  
\item[(iii)] $f_1|_{V}$ is transverse to the subvariety 
$A\subset \widehat Y$ $(V=\C^n\bs Z)$, and
\item[(iv)] the ramification locus of $f_1|_V$ has dimension $\le n-2$.
\end{itemize}

To obtain (iv) we choose $P$ such that $j^1 f_1$ is transverse 
to the subvariety of $J^1(\C^n,Y)$ consisting of all jets of 
non-maximal rank; this subvariety has codimension $n-\dim Y+1\ge 2$
which gives (iv). 

Let $C\subset \C^n$ denote the ramification locus of 
$f_1 \colon \C^n\to \widehat Y$;
thus $\dim C\bs Z \le n-2$ by (iv). The set  
$$
	\Sigma_1 = (\Sigma\cap Z)\cup f_1^{-1}(A) \cup C
$$ 
is an algebraic subvariety of $\C^n$ which does not 
intersect $K$, provided that the approximation of $f$ by $f_1$ is 
sufficiently close on $K$ (we shrink $K$ a little). 
The restriction of $f_1$ to $\C^n\bs \Sigma_1$ 
maps the latter set submersively to 
$Y=\widehat Y\backslash A$.

We claim that $\Sigma_1$ has less $(n-1)$-dimensional 
irreducible components than $\Sigma$. Observe first
that $\dim(\Sigma_1\bs Z)\le n-2$ by properties (iii) and (iv) of $f_1$. 
Next we show $\dim [(Z\bs \Sigma) \cap (f_1^{-1}(A) \cup C)] \le n-2$. 
If $z \in Z \bs \Sigma$ then $f$ is unramified at $z$ 
by the definition of $\Sigma$; furthermore, 
$j^1_{z} f_1 = j^1_{z}f$ by property (ii), 
and hence $f_1$ is also unramified at such point $z$,
thus showing $(Z\bs \Sigma) \cap C =\emptyset$. 
This also implies that $f_1|_{Z\bs \Sigma}$ is transverse to $A$, 
whence $\dim f_1^{-1}(A) \cap (Z\bs \Sigma) \le n-2$.
Hence the $(n-1)$-dimensional 
irreducible components of $Z_1$ are the same as those of 
$\Sigma\cap Z$. Since $z_0\in \Sigma' \backslash Z$, the component 
$\Sigma'$ of $\Sigma$ is not among them which proves the claim.

Repeating the same argument with the pair $(f_1,\Sigma_1)$ gives
an algebraic map $f_2\colon \C^n\to \widehat Y$ and an algebraic
subvariety $\Sigma_2\subset \C^n$ with less $(n-1)$-dimensional
components than $\Sigma_1$ such that 
$f_2\colon \C^n\bs \Sigma_2\to Y$ is a submersion
which approximates $f_1$ (and hence $f$) uniformly on $K$.
Proceeding inductively we obtain in finitely many steps a 
holomorphic submersion 
$\widetilde f\colon \C^n\backslash \wt \Sigma \to Y$, where 
$\wt \Sigma$ is an algebraic subvariety with $\dim \wt \Sigma \le n-2$,
such that $\wt f|_K$ approximates $f|_K$. We complete the proof 
as in Case 1. This establishes Property ${\mathrm S}_{\mathrm n}$ of $Y$.

It remains to prove Property $\mathrm {HS}_{\mathrm n}$ for $n> \dim Y$. 
We shall need the following lemma on algebraic approximation 
of the initial homotopy $f_t$.

\begin{lemma} 
\label{l5.5}
Let $K\subset P\subset \C^n$ be as in $(5.1)$, $(5.2)$.
Let $f_t\colon K\to Y$ for $t\in [0,1]$ be a homotopy of holomorphic 
maps such that $f_0,f_1$ extend to holomorphic maps $P\to Y$.
Let $d$ be a distance function induced by a Riemannian metric on $Y$.
For every $\e>0$ there is an algebraic map 
$F\colon \C^{n+1}\to Y$ such that $d(F(x,t),f_t(x))<\e$
for all $(x,t)\in (K\times [0,1])\cup (P\times \{0,1\}) 
\subset \C^n\times\C$.
\end{lemma}

\begin{proof} 
By corollary \ref{c3.2} we may assume that $f_0$ is algebraic.
Since $Y$ is subelliptic, theorem 4.5 in \cite{FP1} implies 
that we can approximate 
the homotopy $f_t$ uniformly on $K$ by a homotopy 
consisting of holomorphic maps $\widetilde f_t \colon P\to Y$ 
$(t\in [0,1])$ such that $\widetilde f_t=f_t$ for $t=0$ and $t=1$. 
Hence we may assume that the original homotopy $f_t$ 
is defined and holomorphic on $P$, and the initial map 
$f_0$ is algebraic. Lemma \ref{l5.5} now follows from 
theorem \ref{t3.1} applied with the compact set $P$
in $X=\C^n$. 
\end{proof}

Assume that the map $f_t$ in lemma \ref{l5.5} is a submersion 
$K\to Y$ for $t\in [0,1]$, and it is a submersion $P\to Y$ 
for $t=0$ and $t=1$. If the approximation of $f_t$ by 
the algebraic map $F_t=F(\cdotp,t) \colon \C^n\to Y$, 
furnished by lemma \ref{l5.5}, is sufficiently uniformly close on 
the respective sets then we may assume the same properties 
for $F_t$ after a slight shrinking of $K$ and $P$. 

Consider now $F$ as a map to $\widehat Y$. A transversality 
argument, analogous to the one used in the proof of 
Property ${\mathrm S}_{\mathrm n}$, gives a nearby 
algebraic map (still called $F$) such that for all points 
$(x,t)\in \C^{n+1}\bs \Sigma$ outside of an 
algebraic subvariety with $\dim\Sigma\le n-1$ we have 
$F(x,t)\in Y$ and $\partial_x F_t \colon T_x \C^n \to T_{F(x,t)} Y$ 
is surjective.  

Fix such $F$. For all but finitely many values of 
$t\in \C$ the set $\Sigma_t=\{x\in \C^n\colon (x,t)\in\Sigma\}$
then satisfies $\dim \Sigma_t \le n-2$. By a small 
smooth deformation $\tau\colon [0,1]\to\C$ of the parameter interval 
$[0,1] \subset \R\subset \C$ inside $\C$ we can avoid this 
exceptional set of $t$'s and obtain a homotopy of algebraic 
submersions $F_{\tau(t)} \colon \C^n \backslash \Sigma_{\tau(t)} \to Y$ 
with $\dim \Sigma_{\tau(t)} \le n-2$ for all $t\in [0,1]$.
Since $F_{\tau(t)}$ approximates $f_t$ on $K$ (resp.\ on $P$
for $t=0$ and $t=1$), we have $\Sigma_{\tau(t)} \cap K=\emptyset$ for all 
$t\in [0,1]$, and $\Sigma_{\tau(t)}\cap P=\emptyset$ for $t\in \{0,1\}$.
(We shrink $K$ and $P$ slightly.)

Lemma 3.4 in \cite{F3} gives a family of holomorphic 
automorphisms $\psi_t \in\mbox{Aut}\C^n$ depending smoothly on 
$t\in [0,1]$, with $\psi_0$, $\psi_1$ being the identity map, 
such that $\psi_t(P)\cap \Sigma_{\tau(t)} =\emptyset$ 
for every $t\in [0,1]$. The homotopy 
$\widetilde f_t= F_{\tau(t)} \circ \psi_t \colon P\to  Y$ $(t\in [0,1])$
consists of holomorphic submersions approximating 
$f_t$ on $K$ (resp.\ on $P$ for $t\in \{0,1\}$). 
If $\widetilde f_0$ is sufficiently close to $f_0$ 
on $P$ then we can join the two maps by a homotopy of 
submersions $P\to Y$; the same holds at $t=1$. 
\end{proof}

\begin{remark}
Our proof shows that every algebraic manifold $Y$ with a submersive
algebraic spray satisfies Properties $\mathrm{S}_{\mathrm n}$
and $\mathrm {HS}_{\mathrm n}$ (and hence the conclusion
of theorem \ref{t5.1}) for every $n>\dim Y$.
\end{remark}

%
%
%
%
%
\section{Flexibility properties of curves and surfaces}
In this section we survey the holomorphic flexibility properties
of complex curves (Riemann surfaces) and complex surfaces.
Several of the results mentioned here were first proved in
\cite{Gr2}, \cite{F2} and \cite{F6}.

\subsection{Riemann surface} 
We have the following precise result.
The equivalence of (d), (e) and (f) is well known 
and is stated only for completeness.

\begin{proposition} 
\label{p6.1} 
The following are equivalent for a Riemann surface $Y$:
\begin{itemize}
\item[(a)] $Y$ is elliptic (it admits a dominating spray).
\item[(b)] $Y$ enjoys the Oka property.
\item[(c)] Holomorphic submersions $X\to Y$ from any Stein manifold 
$X$ of dimension $\dim X >1$ satisfy the conclusion of theorem \ref{t5.1}.
\item[(d)] $Y$ is dominable by $\C$.
\item[(e)] $Y$ is not hyperbolic.
\item[(f)] $Y$ is one of the surfaces 
$\P_1$, $\C$, $\C^*$ or a torus $\C/\Gamma$.
\end{itemize}
\end{proposition}

\begin{proof} The universal covering of any Riemann surface 
$Y$ is one of the surfaces $\P_1$, $\C$, or 
$\mathbb{U}= \{z\in \C\colon |z|<1\}$. 
$\P_1$ has no nontrivial quotients while $\C$ covers $\C_*$ and 
the tori $\C/\Gamma$. All of these are homogeneous and hence 
admit a dominating spray, and they also satisfy theorem 
\ref{t5.1} since $\P_1$ and $\C$ are of Class $\cA $. 
The disc $\mathbb{U}$ and its quotients are hyperbolic
and hence do not satisfy any flexibility property. 
\end{proof}

%
%
\subsection{Compact complex surfaces}
Let $Y$ be a compact complex surface of Kodaira dimension 
$\kappa$ at most one, the latter condition being implied
by holomorphic dominability of $Y$ by $\C^2$. 
We first consider the case when $Y$ is projective algebraic. 
By the Enriques-Kodaira classification 
(Chapter VI in \cite{BH}) every such $Y$ is obtained from one 
of the following {\it minimal surfaces} $X$ by blowing up points:

(6.2.1) A holomorphic $\P_1$-bundle over a curve 
(Riemann surface) $C$; $\kappa=-\infty$.

(6.2.2) A torus (a quotient of $\C^2$ by a lattice
of real rank four); $\kappa=0$.

(6.2.3) A K3 surface (a compact simply connected surface with 
trivial canonical bundle $K_X=\Lambda^2(T^*X)$); $\kappa=0$.

(6.2.4) A minimal surface $X$ with the structure of an 
{\em elliptic fibration} (\cite{BH}, pp.\ 200-219); 
$\kappa \in \{-\infty,0,1\}$. 

\smallskip

These classes are not entirely disjoint; for example, the moduli 
space of K3 surfaces contains a dense codimension one subset 
consisting of elliptic fibrations. All surfaces on the above list 
are minimal (i.e., not containing any smooth rational $(-1)$-curves), 
and every compact projective surface is obtained from one of them 
(or from a minimal surface of general type) by a finite sequence of blow-ups. 
If ${\rm kod}\, Y\ge 0$ then $Y$ is obtained from a unique minimal $X$. 

A complete list of compact (not necessarily projective algebraic) 
complex surfaces, classified according to the value of 
Kodaira dimension $\kappa<2$, can be found in 
\cite{BH} (Table 10 on p.\ 244). Some of these can be 
included in the classes already listed above;
for example, the class of elliptic fibrations (6.2.4) 
includes all {\em bi-elliptic surfaces} and all 
{\em primary Kodaira surfaces} (\cite{BH}, V.5). 
Besides these one also has

(6.2.5) {\em Secondary Kodaira surfaces}; $\kappa=0$.
These are unramified holomorphic quotients of primary Kodaira surfaces.

(6.2.6) Class VII surfaces ($\kappa(X) =-\infty$, $b_1(X)=1$).
This class includes the {\em Hopf surfaces} (\cite{BH}, V.18),
{\em Inoue surfaces} (\cite{BH}, V.19.), and several others.
\smallskip

%
%
Buzzard and Lu proved in \cite{BL} (theorem 1.1) that
for a compact projective surface $X$ of Kodaira dimension
$\kappa<2$ which is not birational to an elliptic or 
a Kummer K3 surface, dominability by $\C^2$ is equivalent
to the existence of a holomorphic map $\C\to X$ with 
Zariski-dense image. For compact surfaces with $\kappa<2$ 
which are not birationally equivalent to a K3 surface 
they also gave a characterization of dominability in terms 
of the fundamental group (theorem 1.2 in \cite{BL}). 
By the same theorem, every compact projective surface 
which is birationally equivalent to an elliptic or a
Kummer K3 surface is dominable by $\C^2$. The same holds
in the category of compact (non-projective) surfaces 
(propositions 4.4 and 4.5 in \cite{BL}).

In the sequel we consider the question which of the 
compact surfaces from the above list enjoy the Oka property.

{\em Complex tori} (6.2.2), being unramified quotients of $\C^2$,
satisfy the Oka property and the conclusion of theorem \ref{t5.1}.

A {\em secondary Kodaira surface} $X$ (6.2.5) is covered 
by a primary Kodaira surface $\wt X$ which admits the structure of 
an elliptic fibration; hence $X$ enjoys the Oka property if and 
only if $\wt X$ does. There exist surfaces with the structure of a 
ramified elliptic fibration which are not dominable by $\C^2$ 
(example \ref{Winkel1} below).

A {\em Hopf surface} (a special case of a Class VII surface (6.2.6))
is a quotient $\C^2_* /\Gamma$ by a finite extension of a cyclic group. 
It enjoys the Oka property (Corollary 1.5 (i) in \cite{F6})
since its universal covering $\C^2_*=\C^2\bs \{0\}$ is of Class $\cA$.

An {\em Inoue surface} (6.2.6) is a quotient of
$\D\times \C$ (where $\D$ is the unit disc), and hence is 
not dominable by $\C^2$. These surfaces don't admit any 
closed complex curves, and any nonconstant image of $\C$ 
is Zariski dense.

For surfaces (6.2.1) (which include all {\em Hirzebruch surfaces}) 
and for the {\em unramified elliptic fibrations} (6.2.4) 
the following proposition gives a complete answer.

\begin{proposition}
\label{p6.2}
If $C$ is a Riemann surface (not necessarily compact)
and $\pi\colon X\to C$ is either a holomorphic fiber bundle 
with fiber $\P_1$ or an unramified elliptic fibration 
without exceptional fibers then the following are equivalent:
\begin{itemize}
\item[(a)] $X$ satisfies the Oka property.
\item[(b)] $X$ is dominable by $\C^2$.  
\item[(c)] $C$ is not hyperbolic.
\end{itemize}
\end{proposition}

\begin{proof} (a)$\Rightarrow$(b) holds in general (see \S 1), 
(b)$\Rightarrow$(c) is obvious, and (c)$\Rightarrow$(a) 
follows from corollary \ref{c1.6} and proposition \ref{p6.1} 
(in this case $C$ is one of the Riemann surfaces $\P_1$, $\C$, $\C_*$ 
or a torus $\C/\Gamma$). This proves proposition \ref{p6.2}.
We observe in addition that the total space $X$ of a 
holomorphic fiber bundle $X\to C$ with fiber $\P_1$ 
and base $C\in \{\P_1,\C\}$ is of Class $\cA$;
if $C$ is $\C_*$ or a torus $\C/\Gamma$ then $X$ is 
a quotient of a Class $\cA$ manifold. Theorem \ref{t5.1} 
applies to all such manifolds.
\end{proof}

For {\em ramified elliptic fibrations} $\pi\colon X\to C$ 
corollary \ref{c1.6} does not apply, and in general
such $X$ does not enjoy the Oka property even if 
$C$ is non-hyperbolic. The following example was
explained to me by J.\ Winkelmann.

\begin{example}
\label{Winkel1}
There exists a ramified  elliptic fibration $\pi \colon X\to \P_1$
such that $X$ is not dominable by $\C^2$ 
(and hence does not enjoy the Oka property).

Let $Z$ be a hyperelliptic Riemann surface of genus $g\ge 2$,
with involution $\sigma\in {\rm Aut}(Z)$. Then $Z$ is hyperbolic
and $Z/\sigma=\P_1$. Let $E$ be an elliptic curve
considered with its group structure, and let $\tau\in E$
be a non-neutral element of order $2$ ($2\tau=0$ in $E$).
Then $x\to x+\tau$ is a fixed point free involution on $E$.
Take $X=(Z\times E)/\Gamma$ where $\Gamma$ is the cyclic group
of automorphisms generated by $\gamma(z,x)=(\sigma(z), x+\tau)$
which acts without fixed points on $Z\times E$.
Let $\pi\colon X\to \P_1$ be the ramified elliptic fibration 
induced by the projection $p\colon Z\times E\to Z$.  
We have a commutative diagram 
$$ 
	\begin{array}[]{*3c} 
	  Z\times E  & \stackrel{h}{\longrightarrow} & X \\ 
	  \;\;\;\downarrow{p}  & & \;\;\downarrow{\pi} \\
          Z & \longrightarrow & \P_1  
          \end{array}
$$
where $h$ is an unramified holomorphic covering.
A holomorphic map $f\colon \C^n\to X$ lifts to a map
$g\colon \C^n\to Z\times E$ such that $f=h\circ g$. 
Since $Z$ is hyperbolic, the image 
of $g$ is contained in a fiber of $p$, and hence
the map $\pi\circ f\colon \C^n\to \P_1$ is constant.
Thus $X$ is not dominable by $\C^2$. 
(In this example ${\rm kod}\, X=1$.)
\end{example}

This concludes our discussion of elliptic fibrations (6.2.4).

Every compact complex surface bimeromorphic to a 
{\em Kummer K3 surface} (in the class (6.2.3))
is holomorphically dominable by $\C^2$ 
(proposition 4.5 in \cite{BL}). We do not know whether 
any (or all) Kummer K3 surfaces enjoy the Oka property.

\begin{problem}
\label{propermodifications}
To what extent is the Oka property invariant under
proper modifications of (projective) algebraic manifolds? 
In particular, is the Oka property invariant under
blowing up of points and/or blowing down of exceptional
divisors?
\end{problem}

The answer to the first question is affirmative
for manifolds of Class $\cA$ discussed
in \S 2 above, for the simple reason that this class is
closed under blowing up points.

%
%
%
\subsection{Complements of complex subvarieties}
Recall that the complement $\P_p\bs A$ of a thin 
algebraic subvariety is is algebraically subelliptic 
(corollary \ref{c2.4}) and hence enjoys the Oka property.
Call a complex subvariety $A\subset \C^p$ {\em tame} 
if its closure in $\P_p$ does not contain the hyperplane 
at infinity. The complement $\C^p\backslash A$ of a thin tame 
subvariety is elliptic by lemma 7.1 in \cite{FP2}.

Turning to more general surfaces, we have the following
result proved in \cite{F6}. Note that the manifolds $X$ 
in proposition \ref{p6.3} include all complex tori, 
in particular all Abelian varieties.

\begin{proposition} 
\label{p6.3} {\rm (Corollary 1.5 (ii) in \cite{F6})}
Let $X=\C^p/\Gamma$ where $\Gamma$ is a lattice
in $\C^p$ for some $p\ge 2$. 
For any finite set $x_1,\ldots, x_m\in X$ the manifold 
$X_0=X\bs \{x_1,\ldots, x_m\}$ enjoys the Oka property.
\end{proposition}

\begin{proof}
Since \cite{F6} has not been printed yet, we reproduce here the
short proof. Let $\pi\colon \C^p\to X$ denote the quotient projection.
Choose points $q_j\in \C^p$ with $\pi(q_j)=x_j$ for $j=1,\ldots,m$.
The discrete set $\Gamma_0 = \cup_{j=1}^m(\Gamma + q_j)$ 
is tame in $\C^p$ (proposition 4.1 in \cite{BL}, \cite{Bu}) 
and hence $\widetilde X= \C^p\backslash \Gamma_0$ is elliptic.  
Since $\pi\colon \widetilde X \to X_0$ is an unramified 
covering, the Oka property descends from $\wt X$ to $X_0$ 
by corollary \ref{c1.6}.
\end{proof}

The complement of a complex hypersurface of sufficiently 
large degree in $\P_n$ tends to be hyperbolic 
(see e.g.\ \cite{D2}, \cite{D3}, \cite{DE}, 
\cite{DZ}, \cite{Si2}, \cite{SY1}, \cite{SY2}). 
On the other hand, the complement of a smooth cubic curve 
in $\P_2$ is dominable by $\C^2$ according to 
proposition 5.1 in \cite{BL}.

\begin{problem}
\label{complementofcubic}
{\bf (Complements of cubics)}
Does the complement of a smooth cubic curve $C$ in 
$\P_2$ enjoy the Oka property? 
\end{problem}

It can be shown that an affirmative answer 
to the following problem for finitely sheeted maps
$\pi$ would imply an affirmative answer to 
problem \ref{complementofcubic}.

\begin{problem}
\label{descending}
{\bf (Descent of the Oka property)}
Let $\pi\colon X\to X_0$ be a proper holomorphic map
with the ramification locus $\mbox{br}\,\pi$.
Assume that $\pi$ is a subelliptic submersion over 
$X_0\bs \pi(\mbox{br}\,\pi)$ (definition \ref{def1}). 
Does the Oka property of $X$ imply the Oka property of $X_0$?  
What is the answer for finitely sheeted $\pi$? 
\end{problem}

Here is an even more basic question.

\begin{problem}
{\bf (Complements of points)}
Let $X$ be a complex manifold of dimension $\ge 2$
satisfying the Oka property. Does $X\bs \{x_0\}$ 
enjoy the Oka property for every point $x_0\in X$?
\end{problem}

The answer is affirmative for algebraic manifolds
of Class $\cA$ (corollary \ref{c2.4}), 
and the author is not aware of any counterexample.

\medskip
{\it Acknowledgement.}
I express my sincere thanks to G.\ Buzzard, F.\ Kutzschebauch,
F.\ L\'arusson and D.\ Varolin for helpful discussions.
I especially thank J.\ Winkelmann who explained me 
example \ref{Winkel1} and the fact, mentioned in \S 2, 
that an algebraic Lie group without a homomorphism to 
$\C_*$ is algebraically elliptic. I also thank 
D.\ Trotman for pointing out the references \cite{KZ} 
and \cite{MT}, and the referee for thoughtful remarks
which helped me to improve the presentation.

\bibliographystyle{amsplain}

\end{document}